\documentclass[12pt]{article}
\usepackage{authblk}
\usepackage{enumitem}   
\usepackage{enumerate}   
\usepackage{graphicx} 
\usepackage{csquotes}
\usepackage{ragged2e}
\usepackage{varioref}
\usepackage{float}
\usepackage[thinlines]{easytable}
\usepackage{caption}
\usepackage{csquotes}
\usepackage{amsfonts,amssymb,amsmath,exscale,relsize,hyperref,enumitem,amsthm,mathtools,physics, enumerate}
\usepackage[round]{natbib}
\usepackage[usenames]{color}
\usepackage{setspace}
\usepackage{appendix}
\usepackage{graphicx}%
\providecommand{\U}[1]{\protect\rule{.1in}{.1in}}
\newtheorem{theorem}{Theorem}

\newtheorem{proposition}{Proposition}
\newtheorem{remark}{Remark}
\newtheorem{solution*}{Solution}

\def\U{\mathcal{U}}

\def\1{\mathbb{1}}

\newcommand\simiid{\stackrel{iid}{\sim}}
\newcommand\simind{\stackrel{ind}{\sim}}

\providecommand{\keywords}[1]
{
  \small	
  \textbf{\textit{Keywords---}} #1
}

\title{Sharp Asymptotic Minimaxity for Multiple Testing Using One-Group Shrinkage Priors}
\author[1]{Sayantan Paul\thanks{paulsayantan24@gmail.com}}
\author[2]{Prasenjit Ghosh\thanks{prasenjit@stat.tamu.edu}}
\author[3]{Arijit Chakrabarti\thanks{arc@isical.ac.in}}
\affil[1]{Dhirubhai Ambani University, Gandhinagar, Gujarat, India}
\affil[2]{Department of Statistics, Texas A\&M University, College Station, TX, USA}
\affil[3]{Applied Statistics Unit, Indian Statistical Institute, Kolkata, India}

\date{}
\begin{document}

\maketitle

\begin{abstract}

This paper investigates the sharp asymptotic minimaxity properties of Bayesian multiple testing rules in the sparse Gaussian sequence model using a broad class of global--local scale mixtures of normals as priors for the means. Minimaxity is studied under the standard misclassification loss and a composite loss defined as the sum of the false discovery proportion (FDP) and the false non-discovery proportion (FNP). When the sparsity level is known, we show that, by suitably calibrating the global shrinkage parameter to the sparsity level, the proposed testing rule achieves the sharp asymptotic minimax risk under both losses within the \enquote{beta-min} framework. When the sparsity level is unknown, both empirical Bayes and fully Bayesian adaptations of the same testing rule are shown to attain the same sharp asymptotic minimax risk under suitable assumptions on sparsity. Our results reveal that minimaxity is attained for horseshoe-type priors, including the horseshoe, Strawderman--Berger, standard double Pareto, and certain inverse-gamma priors, among others. Additionally, we show that, for priors in the same class with lighter tails, minimaxity fails to hold under either loss function. To the best of our knowledge, these constitute the first sharp asymptotic minimaxity results for multiple testing procedures induced by global--local shrinkage priors.


\end{abstract}

\section{Introduction}
Multiple hypothesis testing has emerged as a critical focal point in statistics, particularly in high-dimensional data analysis. Its relevance stretches across diverse scientific disciplines, including genomics, bioinformatics, medicine, economics, and finance. For example, in microarray experiments, researchers conduct thousands of tests simultaneously to pinpoint differentially expressed genes—those whose expression levels correlate with specific biological traits of interest. This highlights just one canonical example in which analyzing sparse, high-dimensional data is essential, with the primary objective of uncovering significant signals amidst a vast pool of noise. The ongoing exploration of multiple hypothesis testing in the literature not only advances statistical methods but also enhances our ability to extract meaningful insights from complex datasets, underscoring its importance in today's scientific landscape.


Over the past three decades, numerous multiple testing procedures have been proposed in the literature, primarily aimed at controlling an overall measure of type I error at a predetermined level \( \alpha \in (0,1) \). Historically, frequentist multiple testing procedures sought to control the \emph{family-wise error rate} (FWER), defined as the probability of making at least one false rejection. The Bonferroni correction and its refinements ensure FWER control, but are often overly conservative in high-dimensional settings, resulting in very low statistical power (i.e., high type II error). This motivated the search for more flexible error criteria. The seminal work of \cite{BH_1995} introduced the False Discovery Rate (FDR), defined as the expected proportion of false positives among all rejections, as a pragmatic (and less conservative) alternative to the family-wise error rate. The procedure for controlling the FDR, invented by \cite{BH_1995}, is henceforth referred to as the BH procedure. This work was a major impetus for modern research on multiple testing, and a vast body of work soon emerged in both the frequentist and Bayesian frameworks.

While FDR control provides guarantees on type~I errors, it does not directly quantify the overall decision-theoretic performance of a procedure, especially in regimes where signal detection is intrinsically difficult. This has motivated the study of risk criteria that combine some overall type I and type II measures. Prominent examples include the expected misclassification loss (misclassification rate) and composite risk measures such as the expected weighted average of the False Discovery Rate (FDR) and False Non-discovery Rate (FNR). These criteria aim to quantify a more balanced trade-off between false positives and false negatives than approaches that merely control FDR; see, for example, \citet{Genovese_Wasserman_2002}, \citet{Sun_Cai_2007}, \citet{Sun_Cai_2009}, and \citet{BCFG2011}. The overarching goal in this line of work is to identify procedures that achieve favorable global error trade-offs rather than focusing exclusively on one-sided error control.

Parallel to the above lines of development, there has been considerable interest among researchers in the last decade in studying and evaluating hypothesis testing procedures through the prism of minimax analysis for different loss functions.  
While minimax testing has been extensively studied in classical and nonparametric settings, recent work has extended this framework to high-dimensional multiple testing problems \citep{arias2017distribution, belitser2020needles, rabinovich2020optimal, fromont2016family}, raising fundamental questions about optimal detectability under sparsity. Ours is a modest attempt to address some relevant and unresolved questions (to be described in detail later) in this context. 

In minimax analysis, it is essential to consider a sufficiently large parameter space and a broad class of decision rules so that the resulting benchmarks have wide applicability. In hypothesis testing, however, trivial procedures that always reject or never reject yield worst-case risk equal to one under both Hamming loss and FDP$+$FNP loss. Consequently, to obtain meaningful (nontrivial) minimax risks, one must impose conditions that ensure separation between the null and alternative hypotheses so that signals can be detected with positive probability.

In the context of sparse Gaussian sequence models \eqref{eq:1.1}, when the loss function is defined as $\mathrm{FDP}+\mathrm{FNP}$, \citet{arias2017distribution} introduced the \emph{beta-min} condition, requiring non-null effects to exceed a minimum magnitude. Under the additional assumption of equal signal
strengths, they established a sharp dichotomy: below a universal logarithmic threshold, no thresholding rule---including BH-type methods---can achieve nontrivial minimax performance, whereas above this threshold, the BH procedure with a suitably tuned level attains vanishing minimax risk. This highlights the fundamental role of separation conditions in enabling reliable detection. Subsequently, \citet{rabinovich2020optimal} derived non-asymptotic bounds for the combined risk $\mathrm{FDR}+\mathrm{FNR}$ in generalized Gaussian sequence models under additional signal strength assumptions, and showed that BH-type procedures can attain the corresponding lower bounds. However, their signal strength conditions are stronger than the beta-min regime of \citet{arias2017distribution}, leading to asymptotically vanishing risk rather than sharp boundary behavior.

Building on these insights, \citet{abraham2024sharp} fully characterized the phase transition for sparse multiple testing. Under a beta-min assumption and across a wide range of sparsity regimes, they derived explicit asymptotic formulas for the minimax risk under both a composite error criterion (FDP+FNP) and the classification (Hamming) loss, delineating a clean boundary between the impossible and possible regimes of reliable testing. Crucially, they also proved that two popular procedures—BH with a vanishing rejection level and empirical-Bayes $\ell$-value rules (defined in \eqref{eq:lvalue}) under spike–and–slab priors—attain the exact minimax risk \emph{adaptively}, without prior knowledge of the sparsity. Their results extend beyond equal-effect settings, covering heterogeneous signal strengths and more general error models, underscoring the robustness of the boundary phenomenon.

We now turn to Bayesian approaches to multiple testing, which provide a natural decision-theoretic framework through hierarchical modeling and posterior risk minimization. Within this framework, a particularly natural formulation of
simultaneous hypothesis testing is based on two-group mixture models, such as spike-and-slab priors \citep{mitchell1988bayesian}. In the Gaussian sequence model, where the goal is to identify which mean parameters are exactly zero, and which are non-zero, a spike-and-slab prior models each mean as a mixture of a point mass at zero (the \emph{spike}) and a continuous, typically heavy-tailed distribution (the \emph{slab}), with the mixing proportion representing the prior probability of a non-null effect. Under this formulation, the Bayes-optimal decision rule for additive $0$--$1$ loss rejects null hypotheses whose posterior inclusion probabilities are below $1/2$.

This framework has been widely studied and applied in large-scale testing problems; see, for example, \citet{efron2004large}, \citet{Muller_Parmigiani_Robert_Rousseau_2004}, \citet{Scott_Berger_2006}, \citet{Muller_Parmigiani_Rice_2006}, \citet{efron2007large}, and \citet{BCFG2011}, among others. These works establish both the practical effectiveness and the strong theoretical properties of Bayesian two-group models for large-scale multiple testing, including optimality results under sparsity and several risk criteria. More broadly, theoretical properties of two-group mixture formulations for high-dimensional inference have been studied extensively in related contexts; see, for example, \citet{Johnstone_Silverman_2004}, \citet{Cai_Jin_2010}, \citet{Castillo_vanderVaart_2012}, \citet{Neuvial_Roquain_2012}, \citet{Martin_Walker_2014}, and \citet{Castillo_Schmidt_vanderVaart_2015}, only to name a few.

Although intuitive and theoretically appealing on many counts, two-group models often face daunting computational challenges in high-dimensional problems and may be less suitable when most parameters are close to, but not exactly, zero, and posterior sampling or marginal likelihood optimization must be carried out at scale. This has motivated the development of continuous one-group shrinkage priors, which are computationally more tractable and offer greater flexibility. In this formulation, each parameter is modeled as a scale mixture of normals, with a global shrinkage parameter controlling overall sparsity and local shrinkage parameters allowing signal-specific adaptivity. Prominent examples include the inverse Gamma priors \citep{tipping2001sparse}, the horseshoe prior \citep{carvalho2009handling, carvalho2010horseshoe}, the three-parameter beta normal mixtures \citep{Armagan_Clyde_Dunson_2011}, the generalized double Pareto prior \citep{Armagan_Dunson_Lee_2013}, the Dirichlet--Laplace prior \citep{bhattacharya2015dirichlet}, and the horseshoe$+$ prior \citep{bhadra2017horseshoe+}. These priors share the property of \emph{tail robustness}: they strongly shrink small signals toward zero while leaving large effects essentially unshrunk, enabling estimation risks that are within logarithmic factors of the minimax benchmark without prior knowledge of sparsity. See, for instance, \cite{van2014horseshoe}, \cite{van2016conditions}, \cite{van2017adaptive}, \cite{ghosh2017asymptotic}, \cite{bhadra2017horseshoe+}, and \cite{paul2025posterior}, among others. 

In the context of multiple hypothesis testing for the Gaussian sequence model \eqref{eq:1.1}, \citet{carvalho2009handling} proposed an intuitive decision rule based on posterior shrinkage factors arising from the horseshoe prior, which behave analogously to posterior inclusion probabilities in spike-and-slab models. Building on this idea, \citet{datta2013asymptotic} and \citet{ghosh2016asymptotic} studied such rules from a formal decision-theoretic perspective under misclassification (Hamming) loss. They showed that decision rules based on a broad class of heavy-tailed one-group shrinkage priors, including the horseshoe, can achieve asymptotically near-optimal Bayes risk (in the sense of \citealp{BCFG2011}) under sparsity when the data are generated from a two-group model. Subsequent work by \citet{ghosh2017asymptotic} and \citet{paul2025posterior} refined these results and established that testing procedures based on a broad class of so-called \emph{horseshoe-type} priors (defined later) in fact attain exact asymptotic Bayes optimality in terms of type I and type II error trade-offs, paralleling the sharp detection boundaries known for spike-and-slab priors. Collectively, these results demonstrate that carefully chosen one-group shrinkage priors can serve as theoretically optimal, adaptive, and computationally attractive alternatives to discrete mixture models for large-scale multiple testing. Since the two terms misclassification loss and Hamming loss are synonymous, we will use only the term misclassification loss throughout the remaining part of the paper.

To the best of our knowledge, the first investigation of minimax risk properties for
multiple testing procedures based on one-group priors was undertaken by
\citet{salomond2017risk}. That work studied thresholding rules similar to those described
above, using a broad class of Gaussian scale-mixture priors inspired by the analysis of
\citet{van2016conditions}. Under separation conditions analogous to beta-min assumptions
on the non-zero means, \citet{salomond2017risk} derived upper bounds on the maximum risk
measured by the sum of the false discovery rate (FDR) and the false non-discovery rate (FNR).
These results imply that when the non-null signals are sufficiently strong, the
corresponding risk tends to zero asymptotically, indicating consistent detection.

It must be noted that \citet{salomond2017risk} did not derive any expression for asymptotic minimax risk for a class of procedures under suitable conditions on the separation between the null and alternative hypotheses. Consequently, it remains unknown whether the multiple testing rules considered in \citet{salomond2017risk}, or indeed any testing rule based on one-group shrinkage priors, attain the sharp minimax risk established under the beta-min framework of \citet{abraham2024sharp}. Moreover, \citet{salomond2017risk} did not study minimaxity under the classification loss; although an empirical Bayes approach was considered to adapt to unknown sparsity, fully Bayesian procedures were not analyzed. Bridging these gaps between sharp minimax theory and Bayesian one-group testing procedures is the central objective of the present work.

Motivated by the sharp minimax results of \citet{abraham2024sharp}, the unresolved questions in \citet{salomond2017risk}, and the strong Bayes optimality results for one-group priors in \citet{datta2013asymptotic}, \citet{ghosh2016asymptotic}, \citet{ghosh2017asymptotic}, and \citet{paul2025posterior}, we study multiple hypothesis testing using one-group shrinkage priors from the perspective of minimax optimality.

We assume that we have an observation vector, say, $\mathbf{X}=(X_1, \dots, X_n)\in \mathbb{R}^{n}$, from the Gaussian sequence model, such that each $X_i$ can be expressed as
 \begin{equation} \label{eq:1.1}
 	X_i=\theta_i+\epsilon_i, \textrm{ for } i=1,\dots,n. 
 \end{equation}
In \eqref{eq:1.1} above, $\theta_1,\dots,\theta_n$ are unknown mean parameters, and the unobserved residuals $\epsilon_i$'s are assumed to be independent $N(0,\sigma^2)$ random variables. For theoretical developments, the residual variance $\sigma^2$ is typically assumed to be known. Thus, without any loss of generality, one may assume $\sigma^2$ to be $1$. For the above normal means model, we are interested in testing the following hypotheses:
\begin{equation}\label{eq_MHT_problem}
 \textrm{H}_{0i}:\theta_i =0 \textrm{ vs. } \textrm{H}_{1i}:\theta_i \neq 0, \textrm{ simultaneously for } i=1,\ldots,n.   
\end{equation}
Following \citet{abraham2024sharp}, we work under a sparse asymptotic regime in which the parameter vector $\boldsymbol{\theta}=(\theta_{1},\dots,\theta_{n})$ belongs to the set
$$\ell_{0}[q_{n}]=\left\{\theta\in\mathbb{R}^{n}:||\theta||_{0}\leq q_{n} \right\}, \textrm{ with } ||\theta||_{0}:=\#\left\{1\leq i\leq n:\theta_{i}\neq 0\right\},$$
consisting of vectors that have at most $q_{n}$ nonzero coordinates. Throughout the paper, we consider the sparse asymptotic framework where
$$\ q_{n}\rightarrow\infty, \ \textrm{and} \ \dfrac{q_{n}}{n} \rightarrow 0 \textrm{ as } n\rightarrow\infty.$$ We focus on the same broad class of priors previously considered in \cite{ghosh2016asymptotic} and \cite{ghosh2017asymptotic}
when each unknown $\theta_{i}$ is modeled as a scale mixture of normals, as,
\begin{equation} \label{eq:1.6}
\theta_i |\lambda_i,\tau \simind N(0,\lambda^2_i\sigma^2\tau^2), \ \lambda^2_i \simind \pi_1(\lambda^2_i)=K (\lambda^2_i)^{-a-1} L(\lambda^2_i), \ \tau \sim \pi_2(\tau),    
\end{equation}
for some choices of $\pi_2$. In \eqref{eq:1.6} above, $K > 0$ is the constant of proportionality, $a$ is a positive real number, and $L\colon (0,\infty)\mapsto (0,\infty)$ is a measurable, continuous, and slowly varying function. In this formulation, $\lambda_i$'s are called the local shrinkage parameters (calibrating the shrinkage at the individual level) and $\tau$ is called the global shrinkage parameter (inducing the overall sparsity of the parameter vector). As discussed by \citet{polson2010shrink}, heavy-tailed priors on the local scales and substantial prior mass for $\tau$ near zero induce marginal priors on the $\theta_i$'s that resemble two-group mixtures, with strong concentration near zero and heavy tails for large signals.



Considering the composite loss \(\textrm{FDP} + \textrm{FNP}\) and the \emph{hamming} loss, our aim in this paper is to investigate whether the multiple testing rules (formally stated in \eqref{eq:5.chap5inteq16}, \eqref{eq:5.chap5inteq17}, and \eqref{eq:5.chap5part2inteq17}), based on the above class of one-group shrinkage priors, asymptotically achieve the minimax risk under appropriate conditions. One of our theoretical results (Proposition \ref{chap5p2typeIIlower}) indicates that one-group shrinkage priors \eqref{eq:1.6} with \( a > 0.5 \) cannot achieve the desired minimax risk for either loss function. Therefore, we confine our attention to the subclass of \emph{horseshoe-type} priors, where \( a = 0.5 \). This subclass of horseshoe-type priors is quite extensive, encompassing various popular prior distributions such as the horseshoe, Strawderman-Berger, standard double Pareto, and inverse-gamma priors (with shape parameter $0.5$), among others. When the sparsity level is known, we show (Theorems~\ref{chap5p2minmaxtun} and \ref{chap5p2minmaxFDRFNRtun}) that testing rules based on horseshoe-type priors attain the exact minimax risk under both losses, provided the global shrinkage parameter is chosen appropriately as a function of the proportion of non-zero means. When the sparsity level is unknown, we further establish that both empirical Bayes and fully Bayesian adaptations of the same testing rules (Theorems~\ref{chap5p2minmaxEBsimple}, \ref{chap5p2minmaxFB}, \ref{chap5p2minmaxFDRFNREBsimple}, and \ref{chap5p2minmaxFDRFNRFB}) continue to achieve exact minimax risk asymptotically as the number of tests grows to infinity. To the best of our knowledge, these findings are the first of their kind in the literature, demonstrating that properly chosen global-local priors can guarantee exact attainment of minimax risk optimality in the context of a simultaneous hypothesis testing problem.

The remainder of the paper is organized as follows. Section \ref{chap5:sec-exst} reviews existing minimax results in Gaussian sequence models. Section \ref{chap5:sec-ogresl} presents our main theoretical contributions, with Subsection \ref{chap5:sec-ogresl} focusing on classification loss and Subsection \ref{chap5p2minimaxfdrfnr} on composite loss (FDR+FNR). Section \ref{chap5p2conrem}  provides concluding remarks and future directions. Proofs are deferred to Section \ref{chap5p2proofs}.

\subsection{Notations}
\color{black}For any two sequences $\{a_n\}$ and $\{b_n\}$ with $\{b_n\} \neq 0$,
 we write $a_n \asymp b_n$ to denote that there exist two constants $c_1$ and $ c_2$ such that $0< c_1\leq  a_n/b_n \leq c_2<\infty $ for sufficiently large $n$.
Similarly, $a_n \lesssim b_n$ denotes that for sufficiently large $n$, there exists a constant $c>0$ such that $a_n \leq cb_n$. Similarly,
$a_n \propto b_n$ implies there exists some constant $C>0$ such that $a_n= Cb_n$ for all $n$.
\color{black}Finally, $a_n=o(b_n)$ denotes $\lim_{n \to \infty} \frac{a_n}{b_n}=0$. Let $\Phi(\cdot)$ and $\phi(\cdot)$ denote the cumulative distribution function (CDF) and the probability density function (PDF) of the standard normal distribution. Also, $\Bar{\Phi}(x)=1-\Phi(x)$.

\section{Existing Results on Relevant Minimax Risk Theory}
\label{chap5:sec-exst}

Before reviewing the relevant literature on minimax theory in the context of multiple
hypothesis testing, we first introduce some notations and concepts that will be used throughout the rest of this paper. We consider the multiple testing problem \eqref{eq_MHT_problem} under the
Gaussian sequence model \eqref{eq:1.1}.

A multiple testing rule is a measurable function
$\boldsymbol{\psi}=(\psi_1,\ldots,\psi_n)$ taking values in $\{0,1\}^n$, where
$\psi_i=1$ indicates rejection of $H_{0i}$. Given $\boldsymbol{\theta}\in\mathbb{R}^n$
and a testing rule $\boldsymbol{\psi}=(\psi_1,\ldots,\psi_n)$, define the numbers of false discoveries and false non-discoveries as
\[
V(\boldsymbol{\theta},\boldsymbol{\psi})=\sum_{i:\,\theta_i=0}\psi_i, \qquad
T(\boldsymbol{\theta},\boldsymbol{\psi})=\sum_{i:\,\theta_i\neq 0}(1-\psi_i).
\]
Let $R(\boldsymbol{\psi})=\sum_{i=1}^n \psi_i$ denote the total number of rejections and
$S(\boldsymbol{\theta})=\sum_{i=1}^n \mathbf{1}\{\theta_i\neq 0\}$ be the number of true signals.
Then, the false discovery proportion (FDP) and the false non-discovery proportion (FNP) are defined as
\[
\mathrm{FDP}(\boldsymbol{\theta},\boldsymbol{\psi})
=\frac{V(\boldsymbol{\theta},\boldsymbol{\psi})}{\max\{R(\boldsymbol{\psi}),1\}}, \qquad
\mathrm{FNP}(\boldsymbol{\theta},\boldsymbol{\psi})
=\frac{T(\boldsymbol{\theta},\boldsymbol{\psi})}{\max\{S(\boldsymbol{\theta}),1\}}.
\]
 Given $\boldsymbol{\theta}\in\mathbb{R}^n$
and a testing rule $\boldsymbol{\psi}$, the false discovery rate (FDR) is defined as
\begin{align} \label{eq:chap5inteq1}
\mathrm{FDR}(\boldsymbol{\theta},\boldsymbol{\psi})
&= \mathbb{E}_{\boldsymbol{\theta}}\!\left[\mathrm{FDP}(\boldsymbol{\theta},\boldsymbol{\psi})
\right].
\end{align}

Given that a multiple testing rule controls the FDR at a desired level, a natural next
question is whether the same rule can control a measure of type II error. A commonly
used criterion is the false non-discovery rate (FNR), defined as
\begin{align} \label{eq:5.chap5inteq2}
\mathrm{FNR}(\boldsymbol{\theta},\boldsymbol{\psi})
&= \mathbb{E}_{\boldsymbol{\theta}}\!\left[\mathrm{FNP}(\boldsymbol{\theta},\boldsymbol{\psi})\right].
\end{align}

In multiple testing under the sparse regime, sparsity alone does not guarantee reliable detection of true signals. When the nonzero signals are too weak relative to the noise level, the distributions of the observed data under sparse alternatives may be statistically indistinguishable from the distribution under the global null, rendering reliable detection impossible. From a decision-theoretic perspective, this corresponds to a regime in which the minimax risk is no better than that of a trivial data-agnostic rule—specifically, the rule that never rejects any null hypothesis—reflecting insufficient information in the data to distinguish between the null and the alternative. In such regimes, the data are effectively indistinguishable from pure noise, and even the optimal procedure cannot improve upon this trivial benchmark.

To obtain minimax guarantees beyond those achieved by such trivial rules, it is necessary to impose additional separation conditions on the signal strengths beyond sparsity alone so that the alternative is distinguishable from the null. A natural way to impose such separation is through \emph{beta-min} conditions, which require that the nonzero signals exceed a minimum magnitude, thereby ensuring sufficient separation for reliable detection. These conditions delineate the detectability boundary between regimes in which the sparse alternative is statistically indistinguishable from the null and those in which meaningful inference becomes possible. Recent sharp minimax results for the Gaussian sequence model, notably those of \citet{abraham2024sharp}, formalize this principle through explicit beta-min–type conditions calibrated to the ambient dimension and sparsity level.

Given $\tilde a\in\mathbb{R}$, let
\begin{equation}
\boldsymbol{\Theta}(\tilde a,q_n)
= \bigl\{ \boldsymbol{\theta}\in\ell_0[q_n] :
|\theta_i|\ge \tilde a \text{ for all } i\in S_{\boldsymbol{\theta}},\
|S_{\boldsymbol{\theta}}|=q_n \bigr\},
\label{eq:theta_def}
\end{equation}
where $S_{\boldsymbol{\theta}}=\{i:\theta_i\neq 0\}$ denotes the support of
$\boldsymbol{\theta}$. Motivated by earlier works of \citet{arias2017distribution} and
\citet{rabinovich2020optimal}, they chose $\tilde a$ to be close to the universal
threshold $\sqrt{2\log(n/q_n)}$. Specifically, they considered
\begin{equation}
\boldsymbol{\Theta}_b = \boldsymbol{\Theta}(\tilde a_b,q_n),
\qquad
\tilde a_b = \sqrt{2\log\!\left(\frac{n}{q_n}\right)} + b,
\label{eq:theta_b}
\end{equation}
where $b\in\mathbb{R}$ is fixed.

Under this framework, \citet{abraham2024sharp} established that for any $b\in\mathbb{R}$,
\begin{equation}
\inf_{\boldsymbol{\psi}} \sup_{\boldsymbol{\theta}\in\boldsymbol{\Theta}_b}
R(\boldsymbol{\theta},\boldsymbol{\psi})
= 1-\Phi(b) + o(1),
\label{eq:minimax_risk_boundary}
\end{equation}
where
\begin{align} \label{eq:5.chap5inteq3}
R(\boldsymbol{\theta},\boldsymbol{\psi})
= \mathrm{FDR}(\boldsymbol{\theta},\boldsymbol{\psi})
+ \mathrm{FNR}(\boldsymbol{\theta},\boldsymbol{\psi}).
\end{align}
This result remains valid when $b=b_n\to\infty$ or $b=b_n\to -\infty$. It follows that the
oracle rule
\[
\psi_i = \mathbf{1}\!\left\{ |X_i| \ge \sqrt{2\log(n/q_n)} \right\}, \qquad i=1,\ldots,n,
\]
is asymptotically minimax. Moreover, they identified a sharp boundary for nontrivial
testing: when $b\to -\infty$, no procedure outperforms trivial rules, whereas finite $b$
yields nontrivial minimax risk.

Since the oracle rule depends on $q_n$, which is typically unknown, an important
question is whether data-driven procedures can achieve the minimax risk in
\eqref{eq:minimax_risk_boundary}. Assuming $q_n \lesssim n^c$ for some $0<c<1$,
\citet{abraham2024sharp} showed that both the BH procedure with $\alpha_n\to 0$ and
$\ell$-value procedures under spike-and-slab priors with suitable slabs attain the
minimax bound.

Specifically, under the two-group prior
\begin{equation}
\theta_i \stackrel{\mathrm{iid}}{\sim} (1-p)\delta_0 + pF, \ \textrm{ for } i=1,\ldots,n,
\label{eq:two_groups}
\end{equation}
the marginal distribution of $X_i$ is
\begin{equation}
X_i \sim (1-p)\phi(x) + p \int_{\mathbb{R}} \phi(x-\theta)\,F(\theta)\,d\theta,
\label{eq:mixture}
\end{equation}
where $\phi(\cdot)$ denotes the standard normal density. With a quasi-Cauchy slab, they defined
the posterior $\ell$-value
\begin{equation}
\ell_{i,p}(\mathbf{X})
= \Pi_p(\theta_i=0\mid\mathbf{X})
= \frac{(1-p)\varphi(X_i)}{(1-p)\varphi(X_i)+p\tilde f(X_i)},
\label{eq:lvalue}
\end{equation}
and rejected $H_{0i}$ when $\ell_{i,p}(\mathbf{X})<t$. Since $p$ is unknown, they used the
MMLE $\hat p$ and showed that the resulting procedure achieves minimax risk
asymptotically.

\citet{abraham2024sharp} also extended their work to heterogeneous signals. For
$\tilde a=(\tilde a_1,\ldots,\tilde a_{q_n})$, define
\begin{equation}
\boldsymbol{\Theta}(\tilde a,q_n)
= \Bigl\{ \boldsymbol{\theta}\in\ell_0[q_n] :
|\theta_{i_j}|\ge \tilde a_j,\ j=1,\ldots,q_n \Bigr\}.
\label{eq:theta_var}
\end{equation}
They introduced
\[
\Lambda_n(\tilde a)
= \frac{1}{q_n}\sum_{j=1}^{q_n} F_{\tilde a_j}(a_n^*),
\qquad
a_n^*=\bigl(\zeta\log(n/q_n)\bigr)^{1/\zeta},\ \zeta>1,
\]
and proved that for symmetric location models,
\begin{align}
\inf_\psi \sup_{\theta\in\Theta(\tilde a,q_n)}
R(\boldsymbol{\theta},\boldsymbol{\psi})
= \Lambda_n(\tilde a)+o(1).
\label{eq:loc_model}
\end{align}

For the Gaussian sequence model,
$\Lambda_n(\tilde a)=q_n^{-1}\sum_{j=1}^{q_n}[1-\Phi(\tilde a_j-a_n^*)]$ with
$a_n^*=\sqrt{2\log(n/q_n)}$. Both BH and $\ell$-value procedures attain this bound.

They also studied minimax risk under classification (Hamming) loss defined as
\begin{equation}
L(\boldsymbol{\theta},\boldsymbol{\psi})
= \sum_{i=1}^n
\bigl[\mathbf{1}\{\theta_i=0,\psi_i\neq0\}
+ \mathbf{1}\{\theta_i\neq0,\psi_i=0\}\bigr],
\label{eq:hamming}
\end{equation}
and showed that
\begin{equation}
\frac{1}{q_n}\inf_\psi \sup_{\theta\in\Theta_b}
\mathbb{E}L(\boldsymbol{\theta},\boldsymbol{\psi})
= 1-\Phi(b)+o(1),
\label{eq:hamming_minimax}
\end{equation}
which is also achieved by BH and $\ell$-value rules.

As noted earlier, the only prior work addressing minimax risk using one-group priors in the Gaussian sequence model is \citet{salomond2017risk}. Modeling
\begin{equation}
\theta_i\mid\sigma_i^2\sim N(0,\sigma_i^2),
\qquad
\sigma_i^2\sim\pi(\sigma_i^2),
\label{eq:onegroup}
\end{equation}
they proposed thresholding rules of the form $\psi_i = \mathbf{1}\{ \mathbb{E}(\varphi_i \mid \mathbf{X}) > \alpha\}$ with $\varphi_i = \sigma_i^2/(1+\sigma_i^2)$
and obtained upper bounds on FDR+FNR under separation conditions. However, as mentioned before, they did not establish exact minimaxity, did not study misclassification loss, and did not derive sharp minimax boundaries. Thus, the question of minimax optimality of multiple testing rules based on one-group priors remained
open, which we address in the next section.

\section{Theoretical Results}
\label{chap5:sec-ogresl}

In this section, we present the main theoretical results of the paper. As indicated earlier, our interest lies in investigating asymptotic minimaxity properties of multiple testing rules based on one-group priors when the observations arise from the Gaussian sequence model \eqref{eq:1.1}. We begin by formally defining the testing rules. Recall that we model the $\theta_i$'s through the scale-mixture specification in
\eqref{eq:1.6}.

In the hierarchical formulation \eqref{eq:1.6}, the posterior mean of $\theta_i$ admits
the representation
\begin{equation}\label{eq:1.7}
\mathbb{E}(\theta_i\mid \mathbf{X},\tau,\sigma)
=\bigl(1-\mathbb{E}(\kappa_i\mid X_i,\tau,\sigma)\bigr)X_i,
\end{equation}
where $\kappa_i = 1/(1+\lambda_i^2\tau^2)$. Consequently, $\mathbb{E}(1-\kappa_i\mid X_i,\tau,\sigma)$ is the $i$th posterior shrinkage coefficient, quantifying the fraction by which the MLE $X_i$ is shrunk toward zero in the Bayes estimator of $\theta_i$.

By comparing posterior means under spike-and-slab priors and the horseshoe prior, \citet{carvalho2009handling} observed that in sparse regimes, the posterior shrinkage coefficient behaves analogously to the posterior inclusion probability
$\Pi(\theta_i\neq 0\mid \mathbf{X})$ arising from spike-and-slab models with heavy-tailed slabs. This led them to propose the following decision rule (under the horseshoe prior) for the multiple testing problem \eqref{eq_MHT_problem} in model \eqref{eq:1.1}:
\begin{equation}\label{eq:1.17}
\text{Reject } H_{0i}\ \text{ if }\ \mathbb{E}(1-\kappa_i\mid X_i,\tau)>\frac{1}{2},
\ \textrm{ for } i=1,\ldots,n.
\end{equation}
Their simulations yielded promising misclassification performance relative to the Bayes
rule under additive $0$--$1$ loss when the $\theta_i$'s are generated from a spike-and-slab
model. As discussed earlier, this phenomenon was subsequently investigated
theoretically: several authors established formal (near-)optimal asymptotic properties
of related testing rules under a broad class of one-group priors of the form
\eqref{eq:1.6}, when the data are generated from a two-group model. Motivated by these
positive results (albeit not with respect to minimax benchmarks), we investigate the
minimaxity properties of the following testing rule
$\boldsymbol{\psi}=(\psi_1,\ldots,\psi_n)$, where for $i=1,\ldots,n$,
\begin{align}\label{eq:5.chap5inteq16}
\psi_i \equiv \psi_i(X_i)
= \mathbf{1}\Bigl\{\mathbb{E}(1-\kappa_i\mid X_i,\tau)>\frac{1}{2}\Bigr\}.
\end{align}

\paragraph{Assumption on $L(\cdot)$.}
For the theoretical study of the rule \eqref{eq:5.chap5inteq16}, we assume that the
slowly varying function $L(\cdot)$ in \eqref{eq:1.6} satisfies the following condition.

\medskip
\noindent\textbf{\hypertarget{assumption1}{Assumption 1.}}
For $a\ge \frac{1}{2}$:
\begin{itemize}
\item[\textbf{(\hypertarget{A1}{A1})}] There exists $c_0>0$ such that $L(t)\ge c_0$ for all $t\ge t_0$, for some $t_0>0$ depending on $L$ and $c_0$. 
\item[\textbf{(\hypertarget{A2}{A2})}] There exists $M\in(0,\infty)$ such that $\sup_{t\in(0,\infty)} L(t)\le M$.
\end{itemize}
\medskip

Earlier results of \citet{datta2013asymptotic},
\citet{van2014horseshoe}, \citet{ghosh2016asymptotic}, and \citet{ghosh2017asymptotic} clearly indicate that the choice of the global shrinkage parameter $\tau$ plays a pivotal role in capturing the underlying sparsity and in determining asymptotic optimality of the rule
\eqref{eq:5.chap5inteq16} when $\tau$ is tuned using (known) sparsity. In practice, however, the sparsity level is typically unknown. In that case, we replace $\tau$ by an estimate $\widehat{\tau}$ and employ an empirical Bayes version of
\eqref{eq:5.chap5inteq16}:
\begin{align}\label{eq:5.chap5inteq17}
\psi_i \equiv \psi^{\mathrm{EB}}_i(\mathbf{X})
= \mathbf{1}\Bigl\{\mathbb{E}(1-\kappa_i\mid X_i,\widehat{\tau})>\frac{1}{2}\Bigr\},
\ \textrm{ for } i=1,\ldots,n.
\end{align}
Here
\begin{align}\label{eq:5.chap5inteq18}
\widehat{\tau}
= \max\Bigl\{\frac{1}{n},
\frac{1}{c_2 n}\sum_{i=1}^{n}\mathbf{1}\bigl(|X_i|>\sqrt{c_1\log n}\bigr)\Bigr\},
\end{align}
with $c_1\ge 2$ and $c_2\ge 1$, as proposed by \citet{van2014horseshoe}.

An alternative approach under unknown sparsity is to assign an absolutely continuous prior $\pi_2(\tau)$ to $\tau$ (full Bayes). Motivated by \citet{paul2025posterior}, we consider a class of priors $\pi_2(\cdot)$ satisfying:
\medskip

\noindent\textbf{(\hypertarget{C4}{C4})}
\[
\int_{1/n}^{\alpha_n}\pi_2(\tau)\,d\tau=1, \textrm{ for some }\alpha_n \textrm{ such that }\alpha_n\to 0
\ \textrm{and}\ n\alpha_n\to\infty \textrm{ as } n\to\infty.
\]
 \medskip

In the full Bayes procedure, the decision rule is modified as
\begin{align}\label{eq:5.chap5part2inteq17}
\psi_i \equiv \psi^{\mathrm{FB}}_i(\mathbf{X})
= \mathbf{1}\Bigl\{\mathbb{E}(1-\kappa_i\mid \mathbf{X})>\frac{1}{2}\Bigr\},
\ \textrm{ for } i=1,\ldots,n.
\end{align}
where the full Bayes posterior mean of $1-\kappa_i$ is
\begin{align*}
\mathbb{E}(1-\kappa_i\mid \mathbf{X})
&= \int_{1/n}^{\alpha_n}\mathbb{E}(1-\kappa_i\mid \mathbf{X},\tau)\,\pi(\tau\mid \mathbf{X})\,d\tau \\
&= \int_{1/n}^{\alpha_n}\mathbb{E}(1-\kappa_i\mid X_i,\tau)\,\pi(\tau\mid \mathbf{X})\,d\tau .
\end{align*}

We study minimax optimality of the rules \eqref{eq:5.chap5inteq16}, \eqref{eq:5.chap5inteq17}, and \eqref{eq:5.chap5part2inteq17}---depending on whether the sparsity level is known or unknown---under both the classification (Hamming) loss and the composite loss $\mathrm{FDP}+\mathrm{FNP}$. Results under the classification loss are presented in Subsection~\ref{chap5p2minimaxclass}, while Subsection~\ref{chap5p2minimaxfdrfnr} contains results for the $\mathrm{FDP}+\mathrm{FNP}$ loss.

\begin{remark}
    The global parameter $\tau$ controls the overall amount of shrinkage and therefore
plays a crucial role in balancing type~I and type~II errors. If $\tau$ is too small,
the procedure becomes overly conservative and fails to detect signals even at the
beta-min boundary; if $\tau$ is too large, excessive noise coordinates escape
shrinkage, leading to inflated false discoveries. 
In earlier work on Bayes optimality, it has been established that $\tau \asymp q_n/n$ yields positive results. Keeping this in mind, we will consider this as the tuning level of $\tau$ when $q_n$ is known and, in fact, show that it yields an optimal result in this context.
This scaling has also appeared in earlier Bayes-optimality
results for one-group priors and will be adopted throughout our minimax analysis.
\end{remark}
We conclude this preliminary part with a lower bound on the type~II error probability
for the $i$th hypothesis, which explains why the case $a>1/2$ in \eqref{eq:1.6} is not
compatible with minimax optimality for the losses considered here.
\begin{proposition}
\label{chap5p2typeIIlower}
Consider the Gaussian sequence model \eqref{eq:1.1} and the simultaneous hypothesis
testing problem \eqref{eq_MHT_problem} based on the decision rule
\eqref{eq:5.chap5inteq16} using the class of one-group shrinkage priors \eqref{eq:1.6}.
Assume that $\tau\equiv\tau_n\to 0$ as $n\to\infty$ such that
$\frac{n\tau}{q_n}\to C\in(0,\infty)$. Further assume that the slowly varying function
$L(\cdot)$ in \eqref{eq:1.6} satisfies \hyperlink{assumption1}{Assumption~1}. Then for
$a>\frac{1}{2}$ in \eqref{eq:1.6}, for $i=1,\ldots,q_n$, and $b\in\mathbb{R}$,
\begin{align}\label{chap5p2eq101}
\sup_{|\theta_i|\ge \sqrt{2\log(\frac{n}{q_n})}+b} t_{2i} \to 1,
\quad \text{as } n\to\infty,
\end{align}
where $t_{2i}$ denotes the type~II error probability corresponding to the $i$th
hypothesis.
\end{proposition}

Observe that the expected misclassification loss $\mathbb{E}L(\boldsymbol{\theta},\boldsymbol{\psi})$ (see \eqref{eq:hamming}) depends on
both the type~I and type~II error probabilities $t_{1i}$ and $t_{2i}$, respectively. As
will be seen from the proofs of Theorems~\ref{chap5p2minmaxtun}--\ref{chap5p2minmaxFB},
the type~II error is the dominating term in the asymptotic regime considered here.
Hence, the fact that $t_{2i}\to 1$ uniformly over $\boldsymbol{\Theta}_b$ for $a>1/2$
implies that the resulting risk cannot attain the sharp minimax benchmark
$1-\Phi(b)+o(1)$. A similar phenomenon occurs for the composite loss $\mathrm{FDP}+\mathrm{FNP}$: the
proofs of Theorems~\ref{chap5p2minmaxFDRFNRtun}--\ref{chap5p2minmaxFDRFNRFB} show that the
minimax behavior is primarily driven by $\mathrm{FNP}$ (equivalently, by type~II
errors). Consequently, the conclusion of Proposition~\ref{chap5p2typeIIlower} rules out
$a>1/2$ for attaining minimax risk under this loss as well. In what follows, we
therefore focus on the case $a=1/2$ (the \emph{horseshoe-type} class).

\paragraph{On the Case $a < \frac{1}{2}$.}
Proposition~\ref{chap5p2typeIIlower} rules out the case $a>\frac{1}{2}$ by showing that the corresponding testing rules suffer from asymptotically overwhelming type~II errors under the beta-min regime. This naturally raises the question of whether priors with heavier tails, corresponding to $a<\frac{1}{2}$, could also attain the sharp minimax risk.

While we do not pursue a formal impossibility result for this regime in the present paper, several heuristic considerations suggest that substantially different behavior may arise when $a < \frac{1}{2}$. In particular, when $a < \frac{1}{2}$, the prior places substantially larger mass on large local scales, thereby weakening shrinkage of noise coordinates under the null. Consequently, the posterior shrinkage coefficient $\mathbb{E}(1-\kappa_i \mid X_i,\tau)$ may remain non-negligible even for moderate values of $X_i$ arising from null coordinates. This, in turn, may lead to inflated type~I errors near the detection boundary $\sqrt{2\log(n/q_n)}$, thereby obstructing the delicate balance between type~I and type~II errors required for sharp minimax optimality under the beta-min regime.

Motivated by these considerations together with Proposition~\ref{chap5p2typeIIlower}, we restrict attention henceforth to the horseshoe-type class corresponding to $a=\frac{1}{2}$.



\subsection{Results Based On the Classification Loss}
\label{chap5p2minimaxclass}

In this subsection, we establish sharp asymptotic minimaxity under the Classification loss for all three testing rules: the fixed-$\tau$ rule, the EB rule, and the FB rule. We first consider the idealized setting in which the sparsity level is known, and the global shrinkage parameter $\tau$ is tuned accordingly, and then turn to adaptive procedures when sparsity is unknown.
Motivated by Proposition~\ref{chap5p2typeIIlower}, we restrict attention to the case $a=\tfrac12$ in \eqref{eq:1.6}, corresponding to the horseshoe-type class of priors.



\begin{theorem}
\label{chap5p2minmaxtun}
Consider the Gaussian sequence model \eqref{eq:1.1}, and the simultaneous hypothesis testing problem \eqref{eq_MHT_problem} based on the decision rule \eqref{eq:5.chap5inteq16} using the class of one-group shrinkage priors \eqref{eq:1.6}. Assume that $\tau\equiv \tau_{n}\to 0$ as $n\to\infty$ such that $\frac{n\tau}{q_n}\to C\in(0,\infty)$. Further assume that the slowly varying function $L(\cdot)$ in \eqref{eq:1.6} satisfies \hyperlink{assumption1}{Assumption 1}. Then for the sub-class of horseshoe-type priors (i.e.\ $a=\frac{1}{2}$),
\begin{align}\label{eq:chap:5p2.1}
\frac{1}{q_n}\sup_{\boldsymbol{\theta}\in\boldsymbol{\Theta}_b}
\mathbb{E}L(\boldsymbol{\theta},\boldsymbol{\psi})
= 1-\Phi(b)+o(1),
\end{align}
where $L(\boldsymbol{\theta},\boldsymbol{\psi})$ is defined in \eqref{eq:hamming} and $\boldsymbol{\Theta}_b$ is defined in \eqref{eq:theta_b}. This result holds for any $b\in\mathbb{R}$ and also for sequences $b=b_n\to\infty$.
\end{theorem}

\begin{remark}\label{chap5p2minmaxtunrem1}
Theorem~\ref{chap5p2minmaxtun} provides a multiple testing rule based on one-group priors that attains the sharp minimax risk under classification (Hamming) loss. To the best of our knowledge, this is the first such result for global--local priors in this setting. The only condition required on $\tau$ is that it is asymptotically of the same order as $q_n/n$, i.e.\ the proportion of non-null signals. The same scaling was used by \citet{ghosh2016asymptotic} to establish Bayes risk optimality (ABOS-type) when the data are generated from a two-group spike-and-slab model. Our proof proceeds by obtaining appropriate bounds on the type I and type II error probabilities. For the type I error, Theorem~6 of \citet{ghosh2016asymptotic} is useful; the control of the type II error, however, requires additional technical arguments.
\end{remark}

The proof strategy relies on carefully bounding the type~I and type~II error probabilities of the induced testing rule. Control of the type~I error follows from Theorem~6 of \citet{ghosh2016asymptotic}, whereas establishing the corresponding bounds for the type~II error necessitates additional technical arguments tailored to the present minimax framework.

Theorem~\ref{chap5p2minmaxtun} assumes that the sparsity level is known, which is often unrealistic. This motivates the study of the adaptive rules \eqref{eq:5.chap5inteq17} and \eqref{eq:5.chap5part2inteq17}. The next two theorems show that minimax risk under classification loss can still be attained. Their proofs are given in the Section \ref{chap5p2proofs}.

\begin{theorem}
\label{chap5p2minmaxEBsimple}
Consider the setup of Theorem~\ref{chap5p2minmaxtun}, and the empirical Bayes rule \eqref{eq:5.chap5inteq17}, where $\widehat{\tau}$ is defined in \eqref{eq:5.chap5inteq18}. Assume that $q_n \asymp (\log n)^{\delta_2}$ for some $\delta_2>1$. Then for the class horseshoe-type priors (i.e.\ $a=\frac{1}{2}$),
\begin{align}\label{eq:chap:5p2.21}
\frac{1}{q_n}\sup_{\boldsymbol{\theta}\in\boldsymbol{\Theta}_b}
\mathbb{E}L(\boldsymbol{\theta},\boldsymbol{\psi})
= 1-\Phi(b)+o(1),
\end{align}
where $L(\boldsymbol{\theta},\boldsymbol{\psi})$ and $\boldsymbol{\Theta}_b$ are as in Theorem~\ref{chap5p2minmaxtun}. This result holds for any $b\in\mathbb{R}$ and also for sequences $b=b_n\to\infty$.
\end{theorem}

\begin{remark}\label{chap5p2minmaxEBsimplerem1}
Theorem~\ref{chap5p2minmaxEBsimple} answers affirmatively the question of whether the minimax risk under classification loss can still be attained by our decision rule when the sparsity level is unknown. Since $\widehat{\tau}$ depends on the entire dataset, the decision rule for the $i$th hypothesis is no longer a function of $X_i$ alone, which complicates the analysis. Our proof exploits monotonicity properties of $\mathbb{E}(1-\kappa_i\mid X_i,\tau)$ and $\mathbb{E}(\kappa_i\mid X_i,\tau)$ as functions of $\tau$ for fixed $X_i$. For type I error control, we first truncate $\widehat{\tau}$ to a suitable range and apply monotonicity to reduce the problem to bounds derived for fixed $\tau$; outside this range, we use concentration arguments tailored to the structure of $\widehat{\tau}$. For type II error control, we proceed analogously using
the monotonicity of $\mathbb{E}(\kappa_i\mid X_i,\tau)$. Although these ideas resemble those used in Theorems~10--11 of \citet{ghosh2016asymptotic}, the arguments here are necessarily different since we do not assume a two-group data-generating mechanism.
\end{remark}

\begin{remark}\label{chap5p2minmaxEBsimplerem1}
Since $\widehat{\tau}$ depends on the entire dataset, the empirical Bayes decision rule \eqref{eq:5.chap5inteq17} for the $i$th testing problem is no longer a function of $X_i$ alone (as in \eqref{eq:5.chap5inteq16}), which complicates the analysis. The proof of Theorem \ref{chap5p2minmaxEBsimple} exploits monotonicity properties of $\mathbb{E}(1-\kappa_i\mid X_i,\tau)$ and $\mathbb{E}(\kappa_i\mid X_i,\tau)$ as functions of $\tau$, for a fixed $X_i$. For type I error control, we first truncate $\widehat{\tau}$ to a suitable range and use monotonicity to reduce the problem to bounds corresponding to fixed $\tau$; outside this range, we employ concentration arguments tailored to the structure of $\widehat{\tau}$. For type II error control, we proceed analogously, using the monotonicity of $\mathbb{E}(\kappa_i\mid X_i,\tau)$. Although these ideas resemble those in Theorems~10--11 of \citet{ghosh2016asymptotic}, the arguments here are necessarily different, as we do not assume a two-group data-generating mechanism.
\end{remark}

We now present the corresponding optimality result for the full Bayes rule \eqref{eq:5.chap5part2inteq17}.

\begin{theorem}
\label{chap5p2minmaxFB}
Consider the setup of Theorem~\ref{chap5p2minmaxtun}, and the full Bayes rule \eqref{eq:5.chap5part2inteq17}, where the prior on $\tau$ satisfies (\hyperlink{C4}{C4}) with $\alpha_n \propto \frac{(\log n)^{\delta_3}}{n}$ for some $0<\delta_3<\frac{1}{2}$. Assume that $q_n \asymp (\log n)^{\delta_2}$ for some $\delta_2>0$. Then for the sub-class of horseshoe-type priors (i.e.\ $a=\frac{1}{2}$),
\begin{align}\label{eq:chap5p2.21}
\frac{1}{q_n}\sup_{\boldsymbol{\theta}\in\boldsymbol{\Theta}_b}
\mathbb{E}L(\boldsymbol{\theta},\boldsymbol{\psi})
= 1-\Phi(b)+o(1),
\end{align}
where $L(\boldsymbol{\theta},\boldsymbol{\psi})$ and $\boldsymbol{\Theta}_b$ are as in Theorem~\ref{chap5p2minmaxtun}. This result holds for any $b\in\mathbb{R}$ and also for sequences $b=b_n\to\infty$.
\end{theorem}
\begin{remark}
Theorem~\ref{chap5p2minmaxFB} shows that sharp asymptotic minimaxity is preserved even when the global shrinkage parameter is treated fully Bayesianly rather than being fixed in advance or estimated empirically from the data. Thus, exact minimax optimality is not merely a consequence of a carefully calibrated deterministic choice of $\tau$, but remains valid after accounting for uncertainty in $\tau$ through an appropriate prior distribution. In this sense, the sharp minimax behavior established in Theorem~3.1 appears to be an intrinsic property of the horseshoe-type class itself rather than an artefact of a particular tuning strategy. To the best of our knowledge, this is the first sharp asymptotic minimaxity result for a fully Bayesian multiple testing procedure under classification loss based on continuous one-group shrinkage priors.
\end{remark}
\subsection{Results Based on the Composite Loss $\mathrm{FDP}+\mathrm{FNP}$}
\label{chap5p2minimaxfdrfnr}

We now establish sharp asymptotic minimaxity under the composite loss $\mathrm{FDP}+\mathrm{FNP}$ for the same three testing rules. As in the previous subsection, we first analyze the fixed-$\tau$ rule under known sparsity and then show that both the EB and FB rules remain minimax optimal when sparsity is unknown.

\begin{theorem}\label{chap5p2minmaxFDRFNRtun}
Consider the setup of Theorem~\ref{chap5p2minmaxtun}, and the decision rule \eqref{eq:5.chap5inteq16} based on one-group priors for the simultaneous hypothesis testing problem \eqref{eq_MHT_problem}. Assume that $\tau\to 0$ as $n\to\infty$ such that $\frac{n\tau}{q_n}\to C\in(0,\infty)$. Also assume that the slowly varying function $L(\cdot)$ satisfies \hyperlink{assumption1}{Assumption 1}. Then for the subclass of horseshoe-type priors (i.e.\ $a=\frac{1}{2}$),
\begin{align}\label{eq:chap:5p2.51}
\sup_{\boldsymbol{\theta}\in\boldsymbol{\Theta}_b} R(\boldsymbol{\theta},\boldsymbol{\psi})
= 1-\Phi(b)+o(1),
\end{align}
where $\boldsymbol{\Theta}_b$ is as defined in Theorem~\ref{chap5p2minmaxtun} and $R(\boldsymbol{\theta},\boldsymbol{\psi})$ is defined in \eqref{eq:5.chap5inteq3}. This result holds for any $b\in\mathbb{R}$ and also for sequences $b=b_n\to\infty$.
\end{theorem}

\begin{remark}\label{chap5p2minmaxFDRFNRtunrem1}
Theorem~\ref{chap5p2minmaxFDRFNRtun} establishes that the rule \eqref{eq:5.chap5inteq16} also attains the sharp minimax risk for the composite loss $\mathrm{FDP}+\mathrm{FNP}$. In contrast, \citet{salomond2017risk} obtained only upper bounds on the maximum risk for related rules, without establishing exact minimaxity. An advantage of our result is that the minimax risk can be made arbitrarily small by taking $b=b_n\to\infty$. The proof proceeds by separately bounding $\mathrm{FDP}$ and $\mathrm{FNP}$: the truncation used in controlling $\mathrm{FDP}$ is inspired by \citet{salomond2017risk}, but the subsequent bounds require delicate calculations, including the use of Hoeffding's inequality. The bound for $\mathrm{FNP}$ relies on arguments developed for Theorem~\ref{chap5p2minmaxtun} and is independent of \citet{salomond2017risk}.
\end{remark}

Given that \eqref{eq:5.chap5inteq16} attains the minimax risk, it is natural to ask whether the adaptive rules \eqref{eq:5.chap5inteq17} and \eqref{eq:5.chap5part2inteq17} retain minimax optimality when the sparsity level is unknown. The next results answer this question affirmatively. Proofs are provided in the Section \ref{chap5p2proofs}.

\begin{theorem}
\label{chap5p2minmaxFDRFNREBsimple}
Consider the setup of Theorem~\ref{chap5p2minmaxEBsimple}. 
Assume that $q_n \asymp (\log n)^{\delta_2}$ for some $\delta_2>1$, and that $L(\cdot)$ satisfies \hyperlink{assumption1}{Assumption 1}. Then for horseshoe-type priors (i.e.\ $a=\frac{1}{2}$),
\begin{align}\label{eq:chap:5p2.62}
\sup_{\boldsymbol{\theta}\in\boldsymbol{\Theta}_b} R(\boldsymbol{\theta},\boldsymbol{\psi})
= 1-\Phi(b)+o(1),
\end{align}
where $\boldsymbol{\Theta}_b$ is as defined in Theorem~\ref{chap5p2minmaxtun} and $R(\boldsymbol{\theta},\boldsymbol{\psi})$ is defined in \eqref{eq:5.chap5inteq3}. This result holds for any $b\in\mathbb{R}$ and also for sequences $b=b_n\to\infty$.
\end{theorem}

\begin{remark}\label{chap5p2minmaxFDRFNREBsimplerem1}
Theorem~\ref{chap5p2minmaxFDRFNREBsimple} shows that the empirical Bayes rule \eqref{eq:5.chap5inteq17} attains the sharp minimax risk for the normal means model even when the sparsity level is unknown. While \citet{salomond2017risk} proved that an upper bound on the maximum risk converges to $0$ for their proposed rules, our result establishes exact minimaxity within the horseshoe-type subclass. In particular, the minimax risk can be made arbitrarily small by taking $b=b_n\to\infty$.
\end{remark}

Our final result concerns minimax optimality of the full Bayes adaptation \eqref{eq:5.chap5part2inteq17} when the sparsity level is unknown.

\begin{theorem}
\label{chap5p2minmaxFDRFNRFB}
Consider the setup of Theorem~\ref{chap5p2minmaxFB}. Assume that $q_n \asymp (\log n)^{\delta_2}$ for some $\delta_2>0$, and that $L(\cdot)$ satisfies \hyperlink{assumption1}{Assumption 1}. Then for horseshoe-type priors (i.e.\ $a=\frac{1}{2}$),
\begin{align}\label{eq:chap:5p2.620}
\sup_{\boldsymbol{\theta}\in\boldsymbol{\Theta}_b} R(\boldsymbol{\theta},\boldsymbol{\psi})
= 1-\Phi(b)+o(1),
\end{align}
where $\boldsymbol{\Theta}_b$ is as defined in Theorem~\ref{chap5p2minmaxtun} and $R(\boldsymbol{\theta},\boldsymbol{\psi})$ is defined in \eqref{eq:5.chap5inteq3}. This result holds for any $b\in\mathbb{R}$ and also for sequences $b=b_n\to\infty$.
\end{theorem}
\begin{remark}
Theorem~\ref{chap5p2minmaxFDRFNRFB} establishes that the full Bayes rule also attains the sharp minimax risk under the composite loss FDP+FNP. Together with Theorems~\ref{chap5p2minmaxEBsimple}, \ref{chap5p2minmaxFB}, and \ref{chap5p2minmaxFDRFNREBsimple}, this result
shows that exact minimax optimality persists across both empirical Bayes and fully
Bayesian adaptations of the underlying testing rule. Consequently, the sharp minimax
behavior of horseshoe-type priors appears robust to the particular mechanism used to
adapt to unknown sparsity.
\end{remark}

\begin{remark}
The logarithmic sparsity assumptions imposed in Theorems~\ref{chap5p2minmaxEBsimple}, \ref{chap5p2minmaxFB}, \ref{chap5p2minmaxFDRFNREBsimple}, and \ref{chap5p2minmaxFDRFNRFB} arise primarily from the technical arguments used to analyze the adaptive empirical Bayes and fully Bayesian procedures. In particular, the condition $q_n \asymp (\log n)^{\delta_2}$ is inherited from the concentration and posterior-calibration arguments employed in the proofs, rather than from any apparent minimax phenomenon. It would be of interest to determine whether these assumptions can be relaxed while preserving sharp minimax optimality.
\end{remark}

\section{Concluding Remarks}
\label{chap5p2conrem}

This paper provides the first rigorous theoretical demonstration that multiple testing rules induced by a broad class of horseshoe-type global--local shrinkage priors attain the sharp asymptotic minimax risk in the sparse Gaussian sequence model under the beta-min framework of \citet{abraham2024sharp}. We establish this optimality under both the classification (Hamming) loss and the composite loss $\mathrm{FDP}+\mathrm{FNP}$, thereby matching the sharp minimax benchmark previously shown to be attained by the Benjamini--Hochberg procedure and spike-and-slab $\ell$-value rules. These results demonstrate that continuous shrinkage formulations, despite lacking the explicit model-selection mechanism inherent in spike-and-slab approaches, can nevertheless achieve the same sharp minimax limits in sparse multiple testing problems. An additional contribution of this work is clarifying the role of tail behavior within the broad class of one-group shrinkage priors considered in the literature. In particular, Proposition~\ref{chap5p2typeIIlower} shows that priors with $a>1/2$ cannot attain the sharp minimax risk under the beta-min framework, whereas the positive results established in this paper demonstrate that the horseshoe-type subclass attains the exact minimax limits. Consequently, our findings provide new theoretical insight into the role of the tail parameter within this family of priors, a question that remained largely unexplored in earlier optimality studies. Furthermore, we show that these optimality properties persist when the sparsity level is unknown and the global shrinkage parameter is calibrated adaptively, either through an empirical Bayes estimator or through a fully Bayesian prior specification. Thus, the proposed procedures attain the sharp minimax risk without requiring prior knowledge of the underlying sparsity.

The implications of these findings are multifaceted and extend beyond the specific testing rules studied in this paper. At a conceptual level, our results complement a substantial body of earlier work on global--local shrinkage priors. In particular, \citet{van2014horseshoe}, \citet{van2016conditions}, \citet{van2017adaptive}, and \citet{ghosh2017asymptotic} established important theoretical properties of horseshoe-type priors in sparse high-dimensional problems, including posterior contraction, adaptive estimation, asymptotic minimax estimation, and adaptation to unknown sparsity. In the context of multiple testing, \citet{datta2013asymptotic}, \citet{ghosh2016asymptotic}, \citet{ghosh2017asymptotic}, and \citet{paul2025posterior} showed that testing rules induced by one-group shrinkage priors can achieve strong asymptotic Bayes-optimality properties under suitably calibrated two-group data-generating mechanisms. The present work extends this line of research in a fundamentally different direction by establishing sharp asymptotic minimaxity under worst-case parameter configurations within the beta-min framework of \citet{abraham2024sharp}. An interesting consequence is that the same asymptotic calibration $\tau \asymp q_n/n$, which emerged naturally in earlier Bayes-optimality studies, also leads to sharp minimax optimality in the present framework. Taken together, these findings place continuous shrinkage formulations on the same theoretical footing as the Benjamini--Hochberg procedure and spike-and-slab methods from the perspective of sharp asymptotic minimaxity, while further strengthening the theoretical foundations of horseshoe-type priors for high-dimensional inference and multiple testing.

At the same time, several important questions remain open. One natural direction is to investigate whether the marginal maximum likelihood estimator (MMLE) of the global shrinkage parameter $\tau$ \citep{van2017adaptive}, which has played a central role in the estimation literature for global--local priors, also achieves sharp asymptotic minimaxity in the present multiple-testing framework. Another interesting question is whether the minimax theory developed here can be extended beyond the class of horseshoe-type priors to other global--local formulations, including the Dirichlet--Laplace prior, the horseshoe$+$ prior, and related heavy-tailed shrinkage models. It would also be of considerable interest to investigate whether analogous sharp minimax results continue to hold under more general signal configurations and loss functions beyond the beta-min framework considered in this paper. From a broader perspective, extending the present theory beyond the independent Gaussian sequence model to settings involving dependence, generalized linear models, and other non-Gaussian structures remains an important challenge. Finally, the development of fully adaptive procedures that remain minimax optimal when both sparsity and variance parameters are unknown may help further bridge the gap between robust estimation, Bayesian shrinkage methodology, and optimal multiple testing.

In conclusion, our results establish that a broad class of horseshoe-type priors attains the sharp asymptotic minimax risk in sparse multiple testing, thereby placing continuous shrinkage formulations on the same theoretical footing as the Benjamini--Hochberg procedure and spike-and-slab methods. Beyond resolving a fundamental open question concerning the sharp minimaxity of one-group shrinkage priors, the present work strengthens the connection between Bayesian shrinkage methodology and decision-theoretic multiple testing. We hope that these findings will stimulate further research at the intersection of minimax theory, Bayesian adaptivity, and high-dimensional inference.

\section{Proofs}
\label{chap5p2proofs}
\begin{proof}[Proof of Proposition \ref{chap5p2typeIIlower}]
Let us fix any $i \geq 1$. We consider the cases $a \in (\frac{1}{2},1)$ and $a \geq 1$, separately.\\
 
\textbf{Case 1:} Let us first consider the case when $a \in (\frac{1}{2},1)$. Using Theorem 4 of \cite{ghosh2016asymptotic}, and the fact that the slowly varying function $L(\cdot)$ is bounded above by some constant $M \ (> 0)$, we have
\begin{equation}\label{eq:chap5GTGC2016}
  E(1-\kappa_i \big|X_i,\tau) \leq K_{1}e^{\frac{X_i^2}{2}}\tau^{2a}(1+o(1)),
\end{equation}
where the $o(1)$ term above is independent of both the index $i$ and the data point $X_i$, but depends on $\tau$ in such a way that $\lim_{\tau\rightarrow0}o(1)=0$. Now, applying the definition of type II error coupled with \eqref{eq:chap5GTGC2016}, the type II error rate for the $i^{th}$ hypothesis testing problem is given by
    \begin{align} \label{eq:chap5p2typeIIlow1}
        t_{2i} 
        &= \mathbb{P}_{\theta_i \neq 0} \left(\mathbb{E}(1-\kappa_i|X_i,\tau) \leq \frac{1}{2}\right) \nonumber \\
        & \geq \mathbb{P}_{\theta_i \neq 0} \left(K_1 \tau^{2a} \exp (\frac{X^2_i}{2})(1+o(1)) \leq \frac{1}{2}\right) \nonumber \\
        &= \mathbb{P}_{\theta_i \neq 0} \left(X^2_i \leq 4a \log (\frac{1}{\tau})(1+o(1))\right) \nonumber \\
        &= \Phi\left(\sqrt{4a \log (\frac{1}{\tau})}(1+o(1)-\theta_i\right)-\Phi\left(-\sqrt{4a \log (\frac{1}{\tau})}(1+o(1)-\theta_i\right). 
    \end{align}
    
    Note that, to provide a lower bound to the type II error rate on the set $\boldsymbol{\Theta}_b$, it would be enough to consider only one candidate of $\theta_i$ from that set and show that the lower bound of  \eqref{eq:chap5p2typeIIlow1} for that choice of $\theta_i$ tends to $1$ as $n \to \infty$. Let us fix any $b\in\mathbb{R}$ and choose $\theta_i =\sqrt{2 \log (\frac{n}{q_n})}$. Then under the assumption that $\frac{n\tau}{q_n} \to C \in (0,\infty)$ as $n\rightarrow\infty$, we obtain
    \begin{align} 
        \log \left(\frac{1}{\tau}\right)=\log \left(\frac{n}{q_n}\right) + \log \left(\frac{1}{C}\right)+\log (1+o(1)) = \log \left(\frac{n}{q_n}\right) (1+o(1)), 
        \end{align}
       where the $o(1)$ term above is independent of both the index $i$ and the data point $X_i$, but depends on $\tau$ in such a way that $\lim_{\tau\rightarrow0}o(1)=0$. This implies that for all sufficiently large $n$, we have
        \begin{align}
           \sup_{|\theta_i| \geq \sqrt{2 \log (\frac{n}{q_n})}+b} t_{2i}  & \geq \Phi\left((\sqrt{4a}-\sqrt{2})\sqrt{\log \left(\frac{n}{q_n}\right)}\right) (1+o(1))- \Phi\left((-\sqrt{4a}-\sqrt{2})\sqrt{\log \left(\frac{n}{q_n}\right)}\right) (1+o(1)) \nonumber \\
           &= \Phi\left((\sqrt{4a}-\sqrt{2})\sqrt{\log \left(\frac{n}{q_n}\right)}\right) (1+o(1)) \nonumber \\
           & \to 1, \textrm{ as } n \to \infty.
        \end{align}
        In the above chain of inequalities, observe that for any fixed $a \in (\frac{1}{2},1)$, $\Phi\left((\sqrt{4a}-\sqrt{2})\sqrt{\log (\frac{n}{q_n})}\right) (1+o(1)) \to 1$ as $n \to \infty$ and $\Phi\left((-\sqrt{4a}-\sqrt{2})\sqrt{\log (\frac{n}{q_n})}\right) (1+o(1)) \to 0$ as $n \to \infty$.

        \textbf{Case 2:} Next, we consider the case where $a \geq 1$. First, we derive an upper bound to $E(1-\kappa_i|X_i,\tau)$ for any fixed $\tau \in (0,1)$ and $X_i \in \mathbb{R}$.

  Observe that, the posterior distribution of $\kappa_i$ given $X_i$ and $\tau$ is,
	\begin{equation*}
		\pi(\kappa_i|X_i,\tau) \propto \kappa_i^{a-\frac{1}{2}} (1-\kappa_i)^{-a-1}L\left(\frac{1}{\tau^2}(\frac{1}{\kappa_i}-1)\right) \exp{\frac{(1-\kappa)X^2_i}{2}} \hspace*{0.1cm}, \textrm{ for } 0<\kappa_i<1 .
	\end{equation*}
	Using the transformation $t=\frac{1}{\tau^2}(\frac{1}{\kappa}-1)$, $E(1-\kappa_i|X_i,\tau)$ becomes
	\begin{equation*}
		E(1-\kappa_i|X_i,\tau)= \frac{\tau^2 \int_{0}^{\infty} (1+t \tau^2)^{-\frac{3}{2}} t^{-a} L(t) \exp\bigg({\frac{X^2_i}{2}\cdot \frac{t \tau^2}{1+t \tau^2}}\bigg) dt}{\int_{0}^{\infty} (1+t \tau^2)^{-\frac{1}{2}} t^{-a-1} L(t) \exp \bigg({\frac{X^2_i}{2}\cdot \frac{t \tau^2}{1+t \tau^2}} \bigg) dt}.
	\end{equation*}
	Note that, 
	\begin{align*}
		\int_{0}^{\infty} (1+t \tau^2)^{-\frac{1}{2}} t^{-a-1} L(t) \exp\bigg({\frac{X^2_i}{2}\cdot \frac{t \tau^2}{1+t \tau^2}}\bigg) dt & \geq \int_{0}^{\infty} (1+t \tau^2)^{-\frac{1}{2}} t^{-a-1} L(t) dt= K^{-1}(1+o(1)),
	\end{align*}
 where $o(1)$ depends only on $\tau$ and tends to zero as $\tau \to 0$.
	The equality in the last line follows due to the Dominated Convergence Theorem. Hence
	\begin{equation} \label{chap4newL-1.2}
		E(1-\kappa_i|X_i,\tau) \leq K(A_1+A_2)(1+o(1)) \hspace*{0.1cm}, 
	\end{equation}
    where
	\begin{align*}
		A_1 &= \int_{0}^{1} \frac{t \tau^2}{1+t \tau^2} \cdot \frac{1}{\sqrt{1+t \tau^2}} t^{-a-1} L(t) \exp \bigg({\frac{X^2_i}{2}\cdot \frac{t \tau^2}{1+t \tau^2}}\bigg) dt .
	\end{align*}
	and $A_2=  \int_{1}^{\infty} \frac{t \tau^2}{1+t \tau^2} \cdot \frac{1}{\sqrt{1+t \tau^2}} t^{-a-1} L(t) \exp\bigg({\frac{X^2_i}{2}\cdot \frac{t \tau^2}{1+t \tau^2}}\bigg) dt$.
	Since, for any $t \leq 1$ and $\tau \leq 1$, $\frac{t \tau^2}{1+ t \tau^2} \leq \frac{1}{2}$, using the fact $\int_{0}^{\infty} t^{-a-1} L(t) dt=K^{-1}$,
	\begin{equation} \label{chap4newL-1.3}
		A_1 \leq  K^{-1} \tau^2 \exp \bigg({\frac{X^2_i}{4}}\bigg). 
	\end{equation}
	For providing an upper bound on $A_2$, we divide our calculations depending on $a=1$ and $a >1$. \\
    Next, note that, {$1+t \tau^2 \geq t \tau^2$ for any $t>0$ and $t \tau^2 \geq \sqrt{t}$ if and only if $t \geq \frac{1}{\tau^4}$}. As a result, for $a =1$, we have
\begin{align} \label{chap4neweq:5.90}
  &  \int_{1}^{\frac{1}{\tau^4}}  \frac{t \tau^2}{1+t \tau^2} \cdot t^{-2} L(t) dt 
  \leq M \tau^2 \int_{1}^{\frac{1}{\tau^4}} t^{-1} dt= M \tau^2 \log(\frac{1}{\tau^4}), 
\end{align}
     and
     \begin{align}\label{chap4neweq:5.91}
    \int_{\frac{1}{\tau^4}}^{\infty}  \frac{t \tau^2}{1+t \tau^2} \cdot t^{-2} L(t) dt 
  \leq M \tau^2  \int_{\frac{1}{\tau^4}}^{\infty} t^{-\frac{3}{2}} dt =2 M \tau^4.  
\end{align}
Hence, combining \eqref{chap4neweq:5.90} and \eqref{chap4neweq:5.91}, for $a =1$, we obtain
\begin{align}\label{chap4neweq:5.92}
    A_2 & \leq 2 M \exp\bigg(\frac{X^2_i} {2}\bigg) \tau^2[\log(\frac{1}{\tau^4})+\tau^2]= 8M\exp\bigg(\frac{X^2_i} {2}\bigg) \tau^2\log(\frac{1}{\tau})(1+o(1)).
\end{align}
On the other hand, for $a >1$, 
\begin{align}
     \int_{1}^{\infty}  \frac{t \tau^2}{1+t \tau^2} \cdot t^{-a-1} L(t) dt 
  \leq M \tau^2 \int_{1}^{\infty} t^{-a} dt=  \bigg(\frac{M}{(a-1)}\bigg)\tau^2. 
\end{align}
As a result, for $a>1$, 
\begin{align}\label{chap4neweq:5.93}
    A_2 & \leq  \bigg(\frac{M}{(a-1)}\bigg)\tau^2\exp\bigg(\frac{X^2_i} {2}\bigg) .
\end{align}
Therefore, combining \eqref{chap4neweq:5.92} and \eqref{chap4neweq:5.93}, for $a \geq 1$, we have
\begin{align}\label{chap4neweq:5.94}
    A_2 & \leq K_1 \exp\bigg(\frac{X^2_i} {2}\bigg) \tau^2\log(\frac{1}{\tau})(1+o(1)), 
\end{align}
where $K_1$ is a constant depending on $M$ and $a$. Now, combining \eqref{chap4newL-1.2}, \eqref{chap4newL-1.3} and \eqref{chap4neweq:5.94}, we obtain, for $ a \geq 1$,
\begin{align}\label{chap4neweq:5.95}
    E(1-\kappa_i|X_i,\tau) & \leq K_2 \exp\bigg(\frac{X^2_i} {2}\bigg) \tau^2\log(\frac{1}{\tau})(1+o(1)), 
\end{align}
where $K_2$ is a constant depending on $K^{-1}$ and $K_1$. Next, again using the same technique used when $a \in (\frac{1}{2},1)$ with \eqref{chap4neweq:5.95}, a lower bound of $t_{2i}$ is obtained as,
\begin{align} \label{eq:chap4p2typeIIlow1}
        t_{2i} &= \mathbb{P}_{\theta_i \neq 0} (\mathbb{E}(\kappa_i|X_i,\tau) >\frac{1}{2}) \nonumber \\
        & \geq \mathbb{P}_{\theta_i \neq 0} \left( K_2 \exp\bigg(\frac{X^2_i} {2}\bigg) \tau^2\log(\frac{1}{\tau})(1+o(1))\leq \frac{1}{2}\right) \nonumber \\
         &= \mathbb{P}_{\theta_i \neq 0} \left(X^2_i \leq 4 \log (\frac{1}{\tau})-4 \log \log (\frac{1}{\tau})-o(1)-2\log 2\right) \nonumber \\
        &= \mathbb{P}_{\theta_i \neq 0} \left(X^2_i \leq 4 \log (\frac{1}{\tau}) (1+o(1)) \right) \nonumber \\
        &= \Phi\left(\sqrt{4 \log (\frac{1}{\tau})}(1+o(1)-\theta_i\right)-\Phi\left(-\sqrt{4 \log (\frac{1}{\tau})}(1+o(1)-\theta_i\right). 
    \end{align}
    Next, again using the arguments used in the previous case, we have, for all sufficiently large $n$,
        \begin{align}
           \sup_{|\theta_i| \geq \sqrt{2 \log (\frac{n}{q_n})}+b} t_{2i}  & \geq \Phi((\sqrt{4}-\sqrt{2})\sqrt{\log (\frac{n}{q_n})}) (1+o(1))- \Phi((-\sqrt{4}-\sqrt{2})\sqrt{\log (\frac{n}{q_n})}) (1+o(1)) \nonumber \\
           &= \Phi((\sqrt{4}-\sqrt{2})\sqrt{\log (\frac{n}{q_n})}) (1+o(1)) \nonumber \\
           & \to 1, \textrm{ as } n \to \infty.
        \end{align}
\end{proof}

  \begin{proof}[Proof of Theorem \ref{chap5p2minmaxtun}]
        The desired lower bound to the left-hand side of \eqref{eq:chap:5p2.1} follows from Theorem 7 of \cite{abraham2024sharp}. Hence, we only need to obtain an upper bound for it. Using arguments similar to Theorem 7 of \cite{abraham2024sharp}, note that
        \begin{align}  \label{eq:chap:5p2.2}
            \mathbb{E} L(\theta,\psi)= \sum_{i:\theta_i=0} t_{1i}+\sum_{i:\theta_i \neq 0} t_{2i}, 
        \end{align}
        where $t_{1i}$ and $t_{2i}$ denote the type I and type II error rates corresponding to the $i^{th}$ hypothesis testing problem, respectively. Hence, it is enough to show that
        \begin{align} \label{eq:chap:5p2.3}
            \frac{1}{q_n} \sum_{i:\theta_i=0} t_{1i}=o(1),
        \end{align}
        and
        \begin{align}\label{eq:chap:5p2.4}
            \frac{1}{q_n} \sum_{i:\theta_i\neq 0} \sup_{|\theta_i| \geq \sqrt{2 \log (\frac{n}{q_n}) }+b } t_{2i} \leq 1-\Phi(b)+o(1). 
        \end{align}
        First, we prove \eqref{eq:chap:5p2.3}. Using the architecture of the proof of Theorem 6 of \cite{ghosh2016asymptotic} with $a=\frac{1}{2}$, we have
        \begin{align}\label{eq:chap:5p2.5}
            t_{1i} \leq \tilde{K}_1\frac{\tau L\left(\frac{1}{\tau^2}\right)}{\sqrt{\log\left(\frac{1}{\tau^2}\right)}}(1+o(1)), 
        \end{align}
        where $\tilde{K}_1$ is a positive constant that is independent of $\tau$ and $i$, and the $o(1)$ term is such that $\lim_{n \to \infty}o(1)=0$ and is independent of any $i$. Hence, using \eqref{eq:chap:5p2.4} and under the assumption \hyperlink{assumption1}{Assumption 1} on the slowly varying function $L(\cdot)$, we get
        \begin{align}
             \frac{1}{q_n} \sum_{i:\theta_i=0}  t_{1i} & \leq K_1 M \frac{(n-q_n)}{q_n} \frac{\tau}{\sqrt{\log (\frac{1}{\tau^2})}} (1+o(1)) \nonumber \\
             & \leq K_1 M \frac{n \tau}{q_n}  \frac{1}{\sqrt{\log (\frac{1}{\tau^2})}} (1+o(1))\notag\\
            & \rightarrow 0, \textrm{ as } n \to \infty.
        \end{align}
        In the last line of the above chain of inequalities, we used the facts $\frac{n\tau}{q_{n}}\rightarrow C\in(0,\infty)$ as $n\rightarrow\infty$, $\tau\rightarrow 0$ as $n\rightarrow\infty$, and $\log (\frac{1}{\tau^2})\rightarrow\infty$ as $\tau\rightarrow 0$. This completes the proof of \eqref{eq:chap:5p2.3}.\\
        Our next goal is to provide an upper bound to $t_{2i}$ that holds uniformly for each $i=1,\dots,n$. Towards that, we first use Lemma 3 of \cite{ghosh2017asymptotic} that asserts that, for each fixed $x \in \mathbb{R}$ and any fixed $0<\tau<1$, the posterior shrinkage coefficient $\mathbb{E}(\kappa|x,\tau)$ can be bounded above by a non-negative, continuously decreasing, and measurable function $g(x,\tau,\eta,\delta)$ given by 
        $$g(x,\tau,\eta,\delta)=g_1(x,\tau, \eta,\delta)+g_2(x,\tau,\eta,\delta),$$ where 
        $$g_1(x,\tau, \eta,\delta)= C_1\left[x^2 \int_{0}^{\frac{x^2}{1+t_0}} \exp(-\frac{u}{2})du\right]^{-1}=C_1\left[x^2\left\{1-\exp(-\frac{x^2}{2(1+t_0)})\right\}^{-1}\right],$$
        and 
        $$g_2(x,\tau, \eta,\delta)= C_2 \tau^{-1} \exp(-\frac{\eta(1-\delta)}{2}x^2),$$
        for any given $\eta, \delta \in (0,1)$. Here, $C_1$ and $C_2$ are two finite positive constants that do not depend on both $x$ and $\tau$ (but they do depend on $\eta$ and $\delta$), and $t_0$ was previously defined in \hyperlink{assumption1}{Assumption 1}. In addition, the above nonnegative function $g(x,\tau,\eta,\delta)$ satisfies the following:\\
        For any $\rho > \frac{2}{\eta(1-\delta)} $,
        $$\lim_{\tau\rightarrow0}\sup_{|x|>\sqrt{\rho\log(\frac{1}{\tau})}} g(x,\tau,\eta,\delta) = 0.$$
        Let us fix any $\eta \in (0,1)$ and $\delta \in (0,1)$. Choose $\rho =\frac{2}{\eta(1-\delta)}[1+\tilde{g}(\tau)]$ where $\tilde{g}(\tau) >0$ and tends to zero as $\tau \to 0$ such that ${\tau}^{\tilde{g}(\tau)} \to 0$ as $\tau \to 0$. One possible choice of $\tilde{g}(\tau)=(\log (\frac{1}{\tau}))^{-\delta_1}$ for some $0 <\delta_1<1$. Later, we may choose a more specific choice of $\tilde{g}(\tau)$, as required.\\
        Next, we define two events $B_n$ and $C_n$ as $B_n =\{g(X_i, \tau,\eta, \delta) >\frac{1}{2} \}$ and $C_n= \{ |X_i| > \sqrt{\rho \log (\frac{1}{\tau})} \}$ with $\rho$ as defined before. Then, using the preceding observations, the type II error rate for the $i^{th}$ hypothesis testing problem is given by
    \begin{align} \label{eq:chap:5p2.6}
		t_{2i} &= P_{\theta_{i}\neq 0}(E(\kappa_i|X_i,\tau) > \frac{1}{2}) \nonumber\\
		& \leq P_{\theta_{i}\neq 0}(B_n)  \nonumber\\
        & \leq P_{\theta_{i}\neq 0}(B_n\cap C_n)+P_{\theta_{i}}(C^{c}_n). 
	\end{align}
    Observe that, for any $\theta_i\neq 0$
,    \begin{align}  \label{eq:chap:5p2.7}
        P_{\theta_{i}\neq 0}(B_n\cap C_n)  &= P_{\theta_i\neq 0}\left(g(X_i, \tau,\eta, \delta) >\frac{1}{2},|X_i| > \sqrt{\rho \log (\frac{1}{\tau})}\right) \nonumber  \\
        & \leq  P_{\theta_i\neq 0}\left(g_1(X_i,\tau,\eta, \delta) >\frac{1}{4},|X_i| > \sqrt{\rho \log (\frac{1}{\tau})}\right)\nonumber\\ 
        & + P_{\theta_i\neq 0}\left(g_2(X_i, \tau,\eta, \delta) >\frac{1}{4},|X_i| > \sqrt{\rho \log (\frac{1}{\tau})}\right).
    \end{align}
    Since, for given $ \tau,\eta, \delta$, the function $g_{1}(x, \tau,\eta, \delta)$ steadily decreases in $x^2$, for any $|x| >\sqrt{\rho \log (\frac{1}{\tau})}$, we have
    \begin{equation}\notag
        g_{1}(x, \tau,\eta, \delta) \geq \frac{{C}_{1}}{\rho}\left[\log\left(\frac{1}{\tau}\right)\left\{1-\tau^{\frac{\rho}{2(1+t_{0})}}\right\}\right]^{-1} \rightarrow 0 \textrm{ as } \tau\rightarrow 0.
    \end{equation}
  Since $\tau\rightarrow 0$ as $n\rightarrow\infty$, we have, $g_1(X_i, \tau,\eta, \delta)\leq\frac{1}{4}$ with probability $1$ over the set $C_{n}$ for all sufficiently large $n$, and this would be true for any fixed $\eta\in(0,1)$ and $\delta \in (0,1)$, and any $\theta_{i}\neq 0$. Therefore,
  \begin{equation}\label{eq:chap:5p2.7_1}
     P_{\theta_i\neq 0}\left(g_1(X_i, \tau,\eta, \delta)>\frac{1}{4}\mid |X_i| > \sqrt{\rho \log (\frac{1}{\tau})}\right) = 0, \textrm{ for all sufficiently large } n.
  \end{equation}
 Again, we observe that, for given $ \tau,\eta, \delta$, the function $g_{2}(x, \tau,\eta, \delta)$ steadily decreases in $x^2$. Therefore, using a similar argument, it can be shown easily that,
   \begin{equation}\label{eq:chap:5p2.7_2}
     P_{\theta_i\neq 0}\left(g_2(X_i, \tau,\eta, \delta)>\frac{1}{4}\mid |X_i| > \sqrt{\rho \log (\frac{1}{\tau})}\right) = 0, \textrm{ for all sufficiently large } n.
  \end{equation}
 Combining \eqref{eq:chap:5p2.7}, \eqref{eq:chap:5p2.7_1} and \eqref{eq:chap:5p2.7_2}, it therefore follows that 
    \begin{align}\label{eq:chap:5p2.8}
        P_{\theta_{i}\neq 0}(B_n\cap C_n)=0, \textrm{ for all sufficiently large } n.
    \end{align}
 \eqref{eq:chap:5p2.7} and \eqref{eq:chap:5p2.8} together imply that, for all sufficiently large $n$, we have 
\begin{align}\label{eq:chap:5p2.9}
    t_{2i} & \leq P_{\theta_{i}\neq 0}\left(|X_i| \leq \sqrt{\rho \log (\frac{1}{\tau})}\right) \nonumber \\
    &= \Phi(T_n-\theta_i)-\Phi(-T_n-\theta_i), 
\end{align}
 where $T_n = \sqrt{\rho \log (\frac{1}{\tau})} $. Next, we divide our calculations depending on whether $\theta_i$ is positive or not. \\
    \textbf{Case:} Let us assume $\theta_i >0$, that is, $\theta_{i} > \sqrt{2\log\left(\frac{n}{q_{n}}\right)}+b$. Then, using the lower bound on $\theta_i$ coupled with \eqref{eq:chap:5p2.9}, the type II error rate $t_{2i}$ can be bounded above as
    \begin{align} \nonumber
        t_{2i} & \leq \Phi(T_n-\theta_i) \\
        & \leq \Phi\left(\sqrt{\rho \log (\frac{1}{\tau})}-\sqrt{2\log (\frac{n}{q_n})}-b\right),
    \end{align}
    where the choice of $\rho$ is same as discussed before. Observe that the decision rule does not depend on how $\eta \in (0,1), \delta \in (0,1)$ and $\rho >\frac{2}{\eta(1-\delta)}$ is chosen. Hence, taking infimum over all such $\rho$'s and subsequently over all possible choices of $(\eta,\delta) \in (0,1) \times (0,1)$, and finally using the continuity of $\Phi(\cdot)$, we have
    \begin{align} \label{eq:chap:5p2.10}
        t_{2i} & \leq \Phi\left(\sqrt{2(1+g(\tau))\log (\frac{1}{\tau})}-\sqrt{2 \log (\frac{n}{q_n})}-b\right).  
    \end{align}
    Next, under the assumption $\frac{n\tau}{q_n} \to C \in (0,\infty)$, we obtain
    \begin{align} \label{eq:chap:5p2.11}
        \log(\frac{1}{\tau})=\log (\frac{n}{q_n}) + \log (\frac{1}{C})+\log (1+o(1)) 
    \end{align}
    and
    \begin{align} \label{eq:chap:5p2.12}
        g(\tau)\log (\frac{1}{\tau}) = \left(\log (\frac{n}{q_n})\right)^{1-\delta_1}(1+o(1)).  
    \end{align}
    Combining all these observations, we have
    \begin{align} \label{eq:chap:5p2.13}
       t_{2i} & \leq \Phi(\sqrt{a_n+d_n}-\sqrt{a_n}-b)  , 
    \end{align}
    where $a_n= 2\log (\frac{n}{q_n})$ and $d_n=2[\log (\frac{1}{C})+\log(1+o(1))] +2\left(\log(\frac{n}{q_n})\right)^{1-\delta_1}(1+o(1))$. Now, note that
    \begin{align} \label{eq:chap:5p2.14}
        & \Phi(\sqrt{a_n+d_n}-\sqrt{a_n}-b) -\Phi(-b) \nonumber \\
       & = \frac{1}{\sqrt{2 \pi}} \int_{-b}^{\sqrt{a_n+d_n}-\sqrt{a_n}-b} \exp(-\frac{x^2}{2}) dx \nonumber \\
       & \leq  \sqrt{a_n+d_n}-\sqrt{a_n} \nonumber \\
       &= \frac{(\sqrt{a_n+d_n}-\sqrt{a_n})(\sqrt{a_n+d_n}+\sqrt{a_n})} {\sqrt{a_n+d_n}+\sqrt{a_n}} \nonumber \\
       & = \frac{d_n}{\sqrt{a_n+d_n}+\sqrt{a_n}} \nonumber \\
       & \leq \frac{d_n}{2 \sqrt{a_n}}.  
    \end{align}
Note that for $d_n= (\log (\frac{n}{q_n}))^{1-\delta_1}(1+o(1))$, with $\frac{1}{2}<\delta_1<1$, we have
    $d_n= o(\sqrt{a_n})$, as $n \to \infty$. Using this observation, coupled with \eqref{eq:chap:5p2.14}, we obtain
    \begin{align} \label{eq:chap:5p2.15}
        t_{2i} & \leq \Phi(\sqrt{a_n+d_n}-\sqrt{a_n}-b) \nonumber \\
        & = \Phi(-b) + \left[\Phi(\sqrt{a_n+d_n}-\sqrt{a_n}-b)-\Phi(-b)\right] \nonumber \\
        & = \Phi(-b)+o(1). 
    \end{align}
    Note that, this $o(1)$ is independent of $i$.\\
     \textbf{Case 2:} Let us assume that $\theta_i<0$, that is,  $\theta_{i} < -\sqrt{2\log\left(\frac{n}{q_{n}}\right)}-b$. Then, using the lower bound on $-\theta_i$ the upper bound of $t_{2i}$ obtained in \eqref{eq:chap:5p2.9} modifies to
    \begin{align} \label{eq:chap:5p2.16}
        t_{2i} & \leq 1- \Phi(-T_n-\theta_i) \nonumber\\
        & \leq 1- \Phi\left(-\sqrt{\rho \log (\frac{1}{\tau})}+\sqrt{2\log (\frac{n}{q_n})}+b\right). 
    \end{align}
Again using arguments similar to those used in establishing \eqref{eq:chap:5p2.10}-\eqref{eq:chap:5p2.13}, an upper bound to the probability of type II error in this case is obtained as
\begin{align} \label{eq:chap:5p2.17}
       t_{2i} 
       & \leq 1- \Phi(-\sqrt{a_n+d_n}+\sqrt{a_n}+b), 
    \end{align}
    where $a_n$ and $d_n$ are same as before. Next, using the same arguments used to derive \eqref{eq:chap:5p2.14}, we obtain
    \begin{align} \label{eq:chap:5p2.18}
        \Phi(-\sqrt{a_n+d_n}+\sqrt{a_n}+b) = \Phi(b)+o(1), 
    \end{align}
    where $o(1)$ is independent of $i$. Now, combining \eqref{eq:chap:5p2.16}-\eqref{eq:chap:5p2.18}, 
    \begin{align} \label{eq:chap:5p2.19}
        t_{2i} 
       & \leq 1- \Phi(b)+o(1). 
    \end{align}
    Finally \eqref{eq:chap:5p2.15} and \eqref{eq:chap:5p2.19} imply, for each $i$ with $|\theta_i| \geq \sqrt{2 \log (\frac{n}{q_n}) }+b$,
    \begin{align} \label{eq:chap:5p2.20}
        t_{2i} 
       & \leq 1- \Phi(b)+o(1). 
    \end{align}
    This also ensures the upper bound to the left-hand side of \eqref{eq:chap:5p2.4} is of the order of $1- \Phi(b)+o(1)$ and completes the proof of Theorem \ref{chap5p2minmaxtun}.
    \end{proof}

      \begin{proof}[Proof of Theorem \ref{chap5p2minmaxEBsimple}]
       Let $t^{EB}_{1i}$ and $t^{EB}_{2i}$ denote the probabilities of the type I and type II errors for the empirical Bayes decision rule corresponding to the $i^{th}$ hypothesis, respectively. Then the overall risk associated with the empirical Bayes decision rule is given by
        \begin{equation} 
            \mathbb{E} L(\theta,\psi)= \sum_{i:\theta_i=0} t^{EB}_{1i}+\sum_{i:\theta_i \neq 0} t^{EB}_{2i},
        \end{equation}

       Similar to Theorem \ref{chap5p2minmaxtun}, here too, we need to obtain an upper bound to the left hand side of \eqref{eq:chap:5p2.21}. Hence, it suffices to prove that
       \begin{align} \label{eq:chap:5p2.22}
           \sum_{i:\theta_i =0} t^{EB}_{1i}=o(q_n) 
       \end{align}
       and 
       \begin{align}\label{eq:chap:5p2.23}
            \frac{1}{q_n} \sum_{i:\theta_i\neq 0} \sup_{|\theta_i| \geq \sqrt{2 \log (\frac{n}{q_n}) }+b } t^{EB}_{2i} \leq 1-\Phi(b)+o(1). 
        \end{align}
        For any $\alpha_n>0$, we have
        \begin{align}\label{eq:chap:5p2.24}
            t^{EB}_{1i} &=\mathbb{P}_{\theta_i=0}(\mathbb{E}(1-\kappa_i|X_i,\widehat{\tau})>\frac{1}{2}) \nonumber \\
            &= \mathbb{P}_{\theta_i=0}(\mathbb{E}(1-\kappa_i|X_i,\widehat{\tau})>\frac{1}{2}, \widehat{\tau} >2 \alpha_n)+ \mathbb{P}_{\theta_i=0}(\mathbb{E}(1-\kappa_i|X_i,\widehat{\tau})>\frac{1}{2}, \widehat{\tau} \leq 2 \alpha_n) \nonumber \\
            &= u_{1n}+u_{2n}, \textrm{ say, } 
        \end{align}
    where $u_{1n}=\mathbb{P}_{\theta_i=0}(\mathbb{E}(1-\kappa_i|X_i,\widehat{\tau})>\frac{1}{2}, \widehat{\tau} >2 \alpha_n) $ and $u_{2n}=\mathbb{P}_{\theta_i=0}(\mathbb{E}(1-\kappa_i|X_i,\widehat{\tau})>\frac{1}{2}, \widehat{\tau} \leq 2 \alpha_n)$.
      Next, we observe that, for any given $x \in \mathbb{R}$, $\mathbb{E}(1-\kappa|x,\tau)$ is non-decreasing in $\tau$. Hence, using \eqref{eq:chap:5p2.5}, we have
      \begin{align}\label{eq:chap:5p2.25}
          u_{2n} &=\mathbb{P}_{\theta_i=0}(\mathbb{E}(1-\kappa_i|X_i,\widehat{\tau})>\frac{1}{2}, \widehat{\tau} \leq 2 \alpha_n)\nonumber\\
          & = \mathbb{P}_{\theta_i=0}(\mathbb{E}(1-\kappa_i|X_i,\widehat{\tau})>\frac{1}{2}\mid\widehat{\tau} \leq 2 \alpha_n)\mathbb{P}_{\theta_i=0}(\widehat{\tau} \leq 2 \alpha_n) \nonumber\\
          & \leq \mathbb{P}_{\theta_i=0}(\mathbb{E}(1-\kappa_i|X_i,2 \alpha_n)>\frac{1}{2}) \nonumber \\
          & \leq K_1M \frac{2 \alpha_n}{\sqrt{\log (\frac{1}{4 \alpha^2_n})}} (1+o(1)). 
      \end{align}
      Let us define, $\widehat{\tau}_1=\frac{1}{n}$ and $\widehat{\tau}_2=\frac{1}{c_2n}\sum_{i=1}^{n} \mathbf{1}(|X_i| >\sqrt{c_1 \log n})$, with $c_1\geq 2$ and $c_2\geq 1$. Then, $\widehat{\tau}=\max\{\widehat{\tau}_1, \widehat{\tau}_2\}$. Based on this observation, and using the definition of $u_{1n}$, we obtain
      \begin{align}\label{eq:chap:5p2.26}
          u_{1n} & \leq \mathbb{P}_{\theta_i=0}(\widehat{\tau} >2 \alpha_n) \nonumber \\
          & \leq \mathbb{P}_{\theta_i=0}(\widehat{\tau}_1 >2 \alpha_n) +\mathbb{P}_{\theta_i=0}(\widehat{\tau}_2 >2 \alpha_n). 
      \end{align}
      Let us choose $\alpha_n>0$ such that $\widehat{\tau}_1=\dfrac{1}{n}\leq 2 \alpha_n$ for all sufficiently large $n$, so that $\mathbb{P}_{\theta_i=0}(\widehat{\tau}_1 >2 \alpha_n)=0$, provided $n$ is sufficiently large. Hence, \eqref{eq:chap:5p2.26} implies, for all sufficiently large $n$
      \begin{align}\label{eq:chap:5p2.27}
          u_{1n} & \leq \mathbb{P}_{\theta_i=0}(\widehat{\tau}_2 >2 \alpha_n) \nonumber \\
          & \leq \mathbb{P}_{\theta_i=0}(\widehat{\tau}_3 > \alpha_n)+\mathbb{P}_{\theta_i=0}(\widehat{\tau}_4 > \alpha_n), 
      \end{align}
      where $\widehat{\tau}_3= \dfrac{1}{c_2n}\sum\limits_{i: \theta_i \neq 0} \mathbf{1}(|X_i| >\sqrt{c_1 \log n})$ and $\widehat{\tau}_4=\dfrac{1}{c_2n}\sum\limits_{i: \theta_i = 0} \mathbf{1}(|X_i| >\sqrt{c_1 \log n})$. Observe that,
      \begin{align}
          \widehat{\tau}_3 & \leq \frac{1}{c_2} \frac{q_n}{n}.
      \end{align}
      Let us choose $\alpha_n>0$ such that $\alpha_n \geq \frac{1}{c_2}\frac{q_n}{n}$ for all sufficiently large $n$, whence $\mathbb{P}_{\theta_i=0}(\widehat{\tau}_3 > \alpha_n)=0$.

    
      Therefore, for all sufficiently large $n$, we obtain
      \begin{align}\label{eq:chap:5p2.28}
          u_{1n} & \leq \mathbb{P}_{\theta_i=0}(\widehat{\tau}_4 > \alpha_n) \nonumber \\
          &= \mathbb{P}_{\theta_i=0} \left(\frac{1}{n-q_n}\sum_{i: \theta_i = 0} \mathbf{1}(|X_i| >\sqrt{c_1 \log n})> c_2 \frac{n \alpha_n}{n-q_n} \right) \nonumber \\
          & \leq \mathbb{P}_{\theta_i=0}\left(\frac{1}{n-q_n}\sum_{i: \theta_i = 0} \mathbf{1}(|X_i| >\sqrt{c_1 \log n})>  { \alpha_n} \right)\nonumber \\
     &=\mathbb{P}_{\theta_i=0}\left(\frac{1}{n-q_n}\sum_{i: \theta_i = 0} \mathbf{1}(|X_i| >\sqrt{c_1 \log n})> p_{n}+\epsilon_{n} \right),
      \end{align}
 where $p_{n}=2\mathbb{P}_{\theta_i=0}\left(|X_i| >\sqrt{c_1 \log n}\right)$ and $\epsilon_{n}=\alpha_{n}-p_{n}$.\\
 Using Mill's ratio, we obtain,
$$p_{n}=2\mathbb{P}_{\theta_i=0}\left(X_i >\sqrt{c_1 \log n}\right)\propto \frac{n^{-c_{1}/2}}{\sqrt{\log n}}.$$
Again, by our choice of $\alpha_{n}$, we have
\begin{eqnarray}
\epsilon_{n} &\geq& \dfrac{1}{c_{2}}\dfrac{q_{n}}{n}\left(1-\dfrac{\widetilde{C}}{\sqrt{\log n}}\cdot\dfrac{n^{1-c_{1}/2}}{q_{n}}\right) > 0,   
\end{eqnarray}
for all sufficiently large $n$.\\
Using the preceding observations and applying Hoeffding's Inequality, we obtain
      \begin{align}\label{eq:chap:5p2.281}
          u_{1n} & \leq \mathbb{P}_{\theta_i=0}\left(\frac{1}{n-q_n}\sum_{i: \theta_i = 0} \mathbf{1}(|X_i| >\sqrt{c_1 \log n})> p_{n}+\epsilon_{n} \right), \nonumber\\
          & \leq e^{-(n-q_{n})D\left(p_{n}+\epsilon_{n}||p_{n}\right)} \nonumber\\
          & = e^{-(n-q_{n})D\left(\alpha_{n}||p_{n}\right)}, 
      \end{align}
where $D\left(p||q\right)= p\log\left(\dfrac{p}{q}\right)+(1-p)\log\left(\dfrac{1-p}{1-q}\right)$ denotes the Kullback-Leibler divergence between two Bernoulli distributed random variables with parameters $p$ and $q$, respectively.\\
Now, using the facts $\alpha_{n} \geq \frac{1}{c_{2}}\frac{q_{n}}{n}$, and $p_{n}\propto \frac{n^{-c_{1}/2}}{\sqrt{\log n}}$, we obtain

\begin{align}
    \label{chap5minmaxneq3.1}
     D\left(\alpha_{n}||p_{n}\right) &= 
 \alpha_{n}\log\left(\dfrac{\alpha_{n}}{p_{n}}\right)+(1-\alpha_{n})\log\left(\dfrac{1-\alpha_{n}}{1-p_{n}}\right) \nonumber\\
  &= \alpha_{n}\log\left(\dfrac{\alpha_{n}}{p_{n}}\right)-(1-\alpha_{n})\left(\dfrac{\alpha_{n}-p_{n}}{1-p_{n}}\right)(1+o(1)) \nonumber\\
  &\gtrsim \frac{q_{n}}{n}\log(n)(1+o(1)),
\end{align}

whence
\begin{equation}
\label{chap3minmaxeqnew3.3.6}
u_{1n}\leq e^{-(n-q_{n})D\left(\alpha_{n}||p_{n}\right)} \lesssim e^{-(n-q_{n})\frac{q_{n}}{n}\log(n)(1+o(1))} \lesssim e^{-q_{n}\log(n)}(1+o(1)). 
\end{equation}

Let us choose $\alpha_n=c_2\frac{q_n}{n}$, so that $\alpha_n \to 0$ as $n\rightarrow\infty$. The aforesaid choice of $\alpha_{n}$ also satisfies the condition $\alpha_n \geq \frac{1}{c_2}\frac{q_n}{n}$ for all sufficiently large $n$. Hence, combining \eqref{eq:chap:5p2.24},\eqref{eq:chap:5p2.25} and \eqref{chap3minmaxeqnew3.3.6}, we have
\begin{align}\label{eq:chap:5p2.34}
    \frac{1}{q_n} \sum_{i\colon\theta_i=0}  t^{EB}_{1i} & \lesssim \frac{(n-q_n) \alpha_n}{q_n} \frac{1}{\sqrt{\log (\frac{1}{\alpha_n})}} (1+o(1)) + \frac{n}{q_n} e^{-q_n\log(n)}(1+o(1)). 
\end{align}
Note that the choice of $\alpha_n$ satisfies $\frac{n \alpha_n}{q_n} \to C_3$ for some $0<C_3<\infty$. Hence, the first term in the right-hand side of \eqref{eq:chap:5p2.34} goes to zero as $n \to \infty$. On the other hand, since $q_n$ satisfies $ q_n\propto {\log n}^{\delta_2}$ for some $\delta_2>1$, the second term too tends to zero as $n \to \infty$. This implies \eqref{eq:chap:5p2.22} holds. \\
Now we move towards establishing \eqref{eq:chap:5p2.23}. Recall that, by definition, for any given $x \in \mathbb{R}$, $\mathbb{E}(\kappa|x,\tau)$ is non-increasing in $\tau$. Since $\widehat{\tau} \geq \frac{1}{n}$, hence,
\begin{align}\label{eq:chap:5p2.35}
            t^{EB}_{2i} &=\mathbb{P}_{\theta_i\neq 0}(\mathbb{E}(\kappa_i|X_i,\widehat{\tau})>\frac{1}{2}) \nonumber \\
            & \leq  \mathbb{P}_{\theta_i \neq 0}(\mathbb{E}(\kappa_i|X_i,\frac{1}{n})>\frac{1}{2}) 
        \end{align}
Hence, using arguments similar to \eqref{eq:chap:5p2.6}-\eqref{eq:chap:5p2.10} with $\tau=\frac{1}{n}$, when $\theta_i>0$ for $i \in S_{\boldsymbol{\theta}}$, we have,
 \begin{align} \label{eq:chap5p2.10}
        t^{\text{EB}}_{2i} & \leq \Phi(\sqrt{2(1+ (\log n)^{-\delta_1})\log n}-\sqrt{2 \log (\frac{n}{q_n})}-b).  
    \end{align}
 Hence, in order to obtain the desired upper bound on $t^{\text{EB}}_{2i}$, it is sufficient to show that \\ $\sqrt{(1+ (\log n)^{-\delta_1})\log n}-\sqrt{ \log (\frac{n}{q_n})} \to 0$ as $n \to \infty$. Note that, for $q_n \propto (\log n)^{\delta_2}$,
 \begin{align*}
    \sqrt{(1+ (\log n)^{-\delta_1})\log n}-\sqrt{ \log (\frac{n}{q_n})} &=     \sqrt{(1+ (\log n)^{-\delta_1})\log n}- \sqrt{\log n-\delta_2 \log \log n} >0,
 \end{align*}
 for all $n$. Hence, the upper bound converging to zero is enough. Observe that
 \begin{align}\label{eq:chap5p2.101}
      \sqrt{(1+ (\log n)^{-\delta_1})\log n}-\sqrt{ \log (\frac{n}{q_n})}  &= \sqrt{(1+ (\log n)^{-\delta_1})\log n}- \sqrt{\log n-\delta_2 \log \log n} \nonumber \\
      & = \frac{(\log n)^{(1-\delta_1)}+ \delta_2 \log \log n}{\sqrt{(1+ (\log n)^{-\delta_1})\log n} + \sqrt{\log n-\delta_2 \log \log n}} \nonumber \\
      & \leq (\log n)^{(\frac{1}{2}-\delta_1)} + \delta_2 \frac{\log \log n}{\sqrt{\log n}} (1+o(1)) \nonumber \\
      & \to 0, \textrm{ as } n \to \infty, 
 \end{align}
 since $\frac{1}{2} <\delta_1<1$. Combining \eqref{eq:chap5p2.10} and \eqref{eq:chap5p2.101} implies the desired upper bound on $t^{\text{EB}}_{2i}$ is established for $\theta_i >0$ when $i \in S_{\boldsymbol{\theta}}$. The proof when $\theta_i<0$ holds using similar arguments. This completes the proof of Theorem \ref{chap5p2minmaxEBsimple}.
   \end{proof} 

\begin{proof}[Proof of Theorem \ref{chap5p2minmaxFB}]
    Similar to Theorem \ref{chap5p2minmaxEBsimple}, here too, we need to obtain an upper bound of the left-hand side of \eqref{eq:chap5p2.21}. Hence, it suffices to prove that
       \begin{align} \label{eq:chap5p2.22}
           \sum_{i:\theta_i =0} t^{FB}_{1i}=o(q_n) 
       \end{align}
       and 
       \begin{align}\label{eq:chap5p2.23}
            \frac{1}{q_n} \sum_{i:\theta_i\neq 0} \sup_{|\theta_i| \geq \sqrt{2 \log (\frac{n}{q_n}) }+b } t^{FB}_{2i} \leq 1-\Phi(b)+o(1), 
        \end{align}
        where $t^{FB}_{1i}$ and $t^{FB}_{2i}$ are the type I and type II errors for the full Bayes approach corresponding to $i^{th}$ hypothesis, respectively.  Note that,
	\begin{align} \label{eq:T-6.1}
		E(1-\kappa_i|\mathbf{X}) &= \int_{\frac{1}{n}}^{\alpha_n} E(1-\kappa_i|\mathbf{X},\tau) \pi(\tau|\mathbf{X}) d \tau \nonumber \\
		&= \int_{\frac{1}{n}}^{\alpha_n} E(1-\kappa_i|{X}_i,\tau) \pi(\tau|\mathbf{X}) d \tau  \leq E(1-\kappa_i|{X}_i,\alpha_n) \cdot 
	\end{align}
	For proving the inequality above, we first use the fact that given any $X_i \in \mathbb{R}$, $E(1-\kappa_i|X_i,\tau)$ is non-decreasing in $\tau$, and next we use (\hyperlink{C4}{C4}). Now, using Theorem 4 of \cite{ghosh2016asymptotic}, we have for $a=\frac{1}{2}$, 
    as $n \to \infty$
	\begin{equation} \label{eq:T-6.2}
		E(1-\kappa_i|{X}_i,\alpha_n) \leq {K_1} e^{\frac{X^2_i}{2}} {\alpha_n} (1+o(1)).
	\end{equation}
	Here $o(1)$ depends only on $n$ such that $\lim_{n \to \infty}o(1)=0$ and is independent of $i$ and $K_1$ is a constant depending on $M$. Hence, we have the following :
	\begin{align} 
		t^{\text{FB}}_{1i} &= P_{H_{0i}}(E(1-\kappa_i|\mathbf{X})> \frac{1}{2}) \nonumber \\
		& \leq P_{H_{0i}}\bigg(\frac{X^2_i}{2} >  \log (\frac{1}{\alpha_n})- \log K_1- \log (1+o(1)) \bigg) \nonumber \\
		&=  P_{H_{0i}} \bigg(|X_i| > \sqrt{2 \log (\frac{1}{\alpha_n})}\bigg)(1+o(1)) .
	\end{align}
 Note that, as $\alpha_n <1$ for all $n \geq 1$, $2 \log (\frac{1}{\alpha_n}) >0$ for all $n \geq 1$. Next, using the fact that under $H_{0i}, X_i \simiid N(0,1)$ with $1-\Phi(t) < \frac{\phi(t)}{t}$ for $t >0$, we get
 \begin{align} \label{eq:T-6.3}
     t^{\text{FB}}_{1i} &  \leq 2 \frac{\phi(\sqrt{ 2 \log (\frac{1}{\alpha_n})})}{\sqrt{2 \log (\frac{1}{\alpha_n})}}(1+o(1)) = \sqrt{\frac{1}{\pi}} \frac{{\alpha}_n}{\sqrt{\log (\frac{1}{\alpha^2_n}) }} (1+o(1)),
 \end{align}
where $o(1)$ depends only on $n$ such that $\lim_{n \to \infty} o(1)=0$. This implies, for sufficiently large $n$
\begin{align}  \label{eq:T-6.3a}
    \frac{1}{q_n} \sum_{i:\theta_i=0}  t^{\text{FB}}_{1i}  & \lesssim \frac{n \alpha_n}{q_n \sqrt{\log (\frac{1}{\alpha_n}) }} (1+o(1)) \nonumber  \\
    &= \frac{(\log n)^{\delta_3}}{(\log n)^{\delta_2} \sqrt{\log n -\delta_2 \log \log n}}(1+o(1)) \nonumber \\
    &= (\log n)^{\delta_3-\frac{1}{2}-\delta_2}(1+o(1)) \to 0, \textrm{ as } n \to \infty, 
\end{align}
since $\delta_3<\frac{1}{2}$ and $\delta_2>0$. This implies \eqref{eq:chap5p2.22} holds. Next, in order to establish \eqref{eq:chap5p2.23} holds, it is enough to show that, for each $i=1,\cdots,q_n$, $\sup_{|\theta_i| \geq \sqrt{2 \log (\frac{n}{q_n}) }+b } t^{\text{FB}}_{2i} \leq 1-\Phi(b)+o(1)$, where $o(1)$ is independent of $i$. \\
In order to provide an upper bound on the probability of  type-II error induced by the decision rule \eqref{eq:5.chap5part2inteq17}, $t^{\text{FB}}_{2i}= P_{H_{1i}}(E(\kappa_i|\mathbf{X})> \frac{1}{2})$,
 we first note that,
	\begin{align} \label{eq:T-7.4}
		E(\kappa_i|\mathbf{X}) &= \int_{\frac{1}{n}}^{\alpha_n} E(\kappa_i|\mathbf{X},\tau) \pi(\tau|\mathbf{X}) d \tau \nonumber \\
		&= \int_{\frac{1}{n}}^{\alpha_n} E(\kappa_i|{X}_i,\tau) \pi(\tau|\mathbf{X}) d \tau 
		 \leq E(\kappa_i|{X}_i,\frac{1}{n}), 
	\end{align} 
	where the inequality in the last line follows due to the fact that given any $X_i \in \mathbb{R}$, 
 $E(\kappa_i|X_i,\tau)$ is non-increasing in $\tau$. Hence, following the arguments of \eqref{eq:chap5p2.10}-\eqref{eq:chap5p2.101} with $\tau=\frac{1}{n}$, when $\theta_i>0$ for $i=1,\cdots,q_n$, we have,
 \begin{align} \label{eq:chap5p2.110}
      \sup_{|\theta_i| \geq \sqrt{2 \log (\frac{n}{q_n}) }+b }  t^{\text{FB}}_{2i} & \leq 1-\Phi(b)+o(1),  
    \end{align}
    where the $o(1)$ term is independent of any $i$.
This completes the proof of Theorem \ref{chap5p2minmaxFB}.
\end{proof}

   \begin{proof}[Proof of Theorem \ref{chap5p2minmaxFDRFNRtun}]
The lower bound to the left hand side of \eqref{eq:chap:5p2.51} follows from Theorem 1 of \cite{abraham2024sharp}. Hence, it is enough to establish that the upper bound to the left hand side of \eqref{eq:chap:5p2.51} is of the order $1-\Phi(b)+o(1)$. \\
    Using the definition of $FDR(\theta,\psi)$, we have for any $\lambda \in  (0,1)$,
    \begin{align}  \label{eq:chap:5p2.52}
        FDR(\theta,\psi) &= \mathbb{E}_{\theta} \bigg[\frac{\sum_{i:\theta_i=0}\psi_i}{\max(\sum_{i=1}^{n}\psi_i,1)}  \bigg] \nonumber\\
        &= \mathbb{E}_{\theta} \bigg[\frac{\sum_{i:\theta_i=0}\psi_i}{\max(\sum_{i=1}^{n}\psi_i,1)} \mathbf{1}_{\sum_{i=1}^{n}\psi_i > \lambda q_n}\bigg] + \mathbb{E}_{\theta} \bigg[\frac{\sum_{i:\theta_i=0}\psi_i}{\max(\sum_{i=1}^{n}\psi_i,1)} \mathbf{1}_{\sum_{i=1}^{n}\psi_i \leq \lambda q_n}\bigg] \nonumber \\
        &= u_{3n}+u_{4n}, \textrm{ say}. 
    \end{align}
    Observe that, for the set $\mathbf{1}_{\sum_{i=1}^{n}\psi_i > \lambda q_n}$, $\max(\sum_{i=1}^{n}\psi_i,1) > \lambda q_n $. Hence, $u_{3n}$ can be bounded as
    \begin{align} \label{eq:chap:5p2.53}
        u_{3n} & \leq \frac{\sum_{i:\theta_i=0}t_{1i}}{\lambda q_n} = \frac{(n-q_n)t_1}{\lambda q_n} \nonumber \\
        & \leq \frac{M}{\lambda} \frac{n \tau}{q_n}  \frac{(\tau^2)^{a-\frac{1}{2}}}{\sqrt{\log (\frac{1}{\tau^2})}} (1+o(1)) =o(1), \textrm{ as } n \to \infty.
    \end{align}
    The equality in the second line holds since by definition the type I error for $i^{th}$ hypothesis $t_{1i}$ is the same for all $i$, and is denoted as $t_1$. The second inequality follows due to the use of \eqref{eq:chap:5p2.5} and the assumption on $L(\cdot)$.
    Now, note that
    \begin{align} \label{eq:chap:5p2.54}
        \{ \sum_{i=1}^{n} \psi_i \leq \lambda q_n \} \subseteq \{ \sum_{i:\theta_i \neq 0} \psi_i \leq \lambda q_n \}. 
    \end{align}
    This implies
    \begin{align} \label{eq:chap:5p2.55}
        u_{4n} & \leq \mathbb{P}_{\theta}(\sum_{i:\theta_i \neq 0} \psi_i \leq \lambda q_n) .
    \end{align}
    Since $\psi_i$'s are independent Bernoulli random variables, we use Hoeffding's inequality. Choose $0<\lambda <\Phi(b)$. As a result, for all sufficiently large $n$, the constant $t$ defined in Hoeffding's inequality is obtained as
 \begin{align}\label{eq:chap:5p2.56}
     t &= q_n -\sum_{i:\theta_i \neq 0}t_{2i}-\lambda q_n \nonumber \\
     &= (1-\lambda)q_n- \sum_{i:\theta_i \neq 0}t_{2i} \nonumber \\
     & \geq [\Phi(b)-\lambda+o(1)]q_n. 
 \end{align}
 Note that the choice of $\lambda$ implies $t>0$ for sufficiently large $n$. In order to obtain \eqref{eq:chap:5p2.56}, we use a uniform upper bound to all those $t_{2i}$ for which $\theta_i \neq 0$ as derived in \eqref{eq:chap:5p2.20}.
 Using \eqref{eq:chap:5p2.56}, we obtain, for all $|\theta_i| \geq \sqrt{2 \log(\frac{n}{q_n})}+b$, $i \in S_{\boldsymbol{\theta}}$
 \begin{align}\label{eq:chap:5p2.57}
     \mathbb{P}_{\theta}(\sum_{i:\theta_i \neq 0} \psi_i \leq \lambda q_n)  & \leq \exp(-\frac{2 t^2}{q_n}) \nonumber \\
     & \leq \exp(-2[\Phi(b)-\lambda+o(1)]^2q_n). 
 \end{align}
   Next, combining \eqref{eq:chap:5p2.55}-\eqref{eq:chap:5p2.57}, for sufficiently large $n$
   \begin{align}\label{eq:chap:5p2.58}
       \sup_{\theta \in \Theta_b} u_{4n} =o(1), \textrm{ as } n \to \infty. 
   \end{align}
   Now, \eqref{eq:chap:5p2.52}, \eqref{eq:chap:5p2.53} and \eqref{eq:chap:5p2.58} combine to prove that
   \begin{align}\label{eq:chap:5p2.59}
       \sup_{\theta \in \Theta_b} FDR(\theta,\psi) =o(1), \textrm{ as } n \to \infty. 
   \end{align}
   On the other hand, for $FNR(\theta,\psi)$, from definition, we obtain
   \begin{align}\label{eq:chap:5p2.60}
       FNR(\theta,\psi) &= \frac{\mathbb{E}_{\theta}(\sum_{i:\theta_i \neq 0}(1-\psi_i))}{q_n} \nonumber \\
       &= \frac{1}{q_n} \sum_{i:\theta_i \neq 0} t_{2i}. 
   \end{align}
   Recall that, in \eqref{eq:chap:5p2.20}, we have already derived that, for each $|\theta_i| \geq \sqrt{2 \log (\frac{n}{q_n}) }+b$, $i \in S_{\boldsymbol{\theta}}$,
    \begin{align}
        t_{2i} 
       & \leq 1- \Phi(b)+o(1),
    \end{align}
    where the $o(1)$ term is independent of any $i$.
    This, in combination with \eqref{eq:chap:5p2.60}, implies that, for all sufficiently large $n$,
    \begin{align} \label{eq:chap:5p2.61}
         \sup_{\theta \in \Theta_b} FNR(\theta,\psi) & \leq 1- \Phi(b)+o(1). 
    \end{align}
    Finally, \eqref{eq:chap:5p2.59} and \eqref{eq:chap:5p2.61} together imply the upper bound to the left-hand side of \eqref{eq:chap:5p2.51} is of the order $1- \Phi(b)+o(1)$ and completes the proof.
\end{proof}
\begin{proof}[Proof of Theorem \ref{chap5p2minmaxFDRFNREBsimple}]
   Here we follow the basic architecture of the proof of Theorem \ref{chap5p2minmaxFDRFNRtun}. Hence, we only need to show that the upper bound to the left-hand side of \eqref{eq:chap:5p2.62} is of the order $1-\Phi(b)+o(1)$. \\
   Next, using the argument of \eqref{eq:chap:5p2.52}, we have
   \[
u_{5n}
=
\mathbb{E}_{\boldsymbol{\theta}}
\left[
\frac{\sum_{i:\theta_i=0}\psi_{i}^{\textrm{EB}}}{\max\left(\sum_{i=1}^n \psi_{i}^{\textrm{EB}},1\right)}
\mathbf{1}\left\{\sum_{i=1}^n \psi_{i}^{\textrm{EB}}>\lambda q_n\right\}
\right],
\]
\[
u_{6n}
=
\mathbb{E}_{\boldsymbol{\theta}}
\left[
\frac{\sum_{i:\theta_i=0}\psi_{i}^{\textrm{EB}}}{\max\left(\sum_{i=1}^n \psi_{i}^{\textrm{EB}},1\right)}
\mathbf{1}\left\{\sum_{i=1}^n \psi_{i}^{\textrm{EB}}\le \lambda q_n\right\}
\right],
\]
and $\psi_{i}^{\textrm{EB}}$ is defined in (21).
   \begin{align} \label{eq:chap:5p2.63}
        u_{5n} & \leq \frac{\sum_{i:\theta_i=0}t^{EB}_{1i}}{\lambda q_n} \nonumber  \\
        & \lesssim \frac{n \alpha_n}{q_n} \frac{1}{\sqrt{\log (\frac{1}{\alpha_n})}} (1+o(1)) + \frac{n}{q_n} e^{-\frac{ q_n}{c_2}(1+o(1))}. 
    \end{align}
   Here inequality in the second line holds due to \eqref{eq:chap:5p2.34}. Note that the choice of $\alpha_n$ satisfies $\frac{n \alpha_n}{q_n} \to C_3$ for some $0<C_3<\infty$. Hence, the first term to the right hand side of \eqref{eq:chap:5p2.63} goes to zero as $n \to \infty$. On the other hand, since $q_n$ satisfies $ q_n\propto \left({\log n}\right)^{\delta_2}$ for some $\delta_2>1$, the second term too tends to zero as $n \to \infty$. This implies
   \begin{align}\label{eq:chap:5p2.64}
       u_{5n}=o(1), \textrm{ as } n \to \infty. 
   \end{align}
Since the ratio term in the definition of \(u_{6n}\) is bounded above by one,
\[
u_{6n}
\le
\mathbb{P}_{\boldsymbol{\theta}}\left(
\sum_{i=1}^n \psi_i^{\textrm{EB}}\le \lambda q_n
\right)
\le
\mathbb{P}_{\boldsymbol{\theta}}\left(
\sum_{i:\theta_i\ne0}\psi_i^{\textrm{EB}}\le \lambda q_n
\right).
\]

Let
\[
\psi_i^{(0)}
=
\mathbf{1}\left\{
\mathbb{E}(1-\kappa_i\mid X_i,1/n)>\frac12
\right\}.
\]
Since $\widehat\tau\ge 1/n$ and $\mathbb{E}(1-\kappa_i\mid X_i,\tau)$ is non-decreasing in $\tau$, we have
\[
\psi_{i}^{\textrm{EB}}\ge \psi_i^{(0)}, \qquad i=1,\ldots,n.
\]
Hence,
\[
\left\{\sum_{i:\theta_i\ne0}\psi_{i}^{\textrm{EB}}\le \lambda q_n\right\}
\subseteq
\left\{\sum_{i:\theta_i\ne0}\psi_i^{(0)}\le \lambda q_n\right\}.
\]
Let
\[
S^{(0)}=\sum_{i:\theta_i\ne0}\psi_i^{(0)}
\]
and denote by $t_{2i}^{(0)}$ the type~II error probability corresponding to the fixed rule $\psi_i^{(0)}$. Then
\[
\mathbb{E}_{\boldsymbol{\theta}}[S^{(0)}]
=
\sum_{i:\theta_i\ne0}(1-t_{2i}^{(0)}).
\]
The variables $\psi_i^{(0)}$, $i\in S_{\boldsymbol{\theta}}$, are independent Bernoulli random variables. Therefore, applying Hoeffding's inequality, we obtain
\[
u_{6n}
\le
\mathbb{P}_{\boldsymbol{\theta}}\left(S^{(0)}\le \lambda q_n\right)
=
\mathbb{P}_{\boldsymbol{\theta}}\left(S^{(0)}\le \mathbb{E}_{\boldsymbol{\theta}}[S^{(0)}]-t\right),
\]
where
\begin{align}
t
=
q_n - \sum_{i:\theta_i\ne 0} t^{(0)}_{2i} - \lambda q_n.
\label{eq:chap:5p2.571}
\end{align}
Using the same type~II error bound as in empirical Bayes rule, applied to the fixed rule $\psi_i^{(0)}$, we obtain, for sufficiently large $n$,
\begin{align}
t
&=
q_n-\sum_{i:\theta_i\ne 0} t^{(0)}_{2i}-\lambda q_n \nonumber\\
&\ge
[\Phi(b)-\lambda+o(1)]q_n.
\label{eq:chap:5p2.572}
\end{align}
   Again employing \eqref{eq:chap:5p2.54} and \eqref{eq:chap:5p2.55}, we obtain for sufficiently large $n$
   \begin{align}
        t &= q_n -\sum_{i:\theta_i \neq 0}t^{(0)}_{2i}-\lambda q_n \nonumber \\
     & \geq [\Phi(b)-\lambda+o(1)]q_n, 
   \end{align}
   where the upper bound on $t^{EB}_{2i}$ as obtained in \eqref{eq:chap5p2.10} is used.
The choice $0<\lambda<\Phi(b)$ ensures that the right-hand side is positive for all sufficiently large $n$. Hence, Hoeffding's inequality gives for sufficiently large $n$
   \begin{align}\label{eq:chap:5p2.66}
       \sup_{\theta \in \Theta_b} u_{6n} =o(1), \textrm{ as } n \to \infty. 
   \end{align}
   Combining \eqref{eq:chap:5p2.64} and \eqref{eq:chap:5p2.66} ensures
   \begin{align}\label{eq:chap:5p2.67}
       \sup_{\theta \in \Theta_b} FDR(\theta,\psi) =o(1), \textrm{ as } n \to \infty. 
   \end{align}
   For $FNR(\theta,\psi)$, from definition, we obtain
   \begin{align}\label{eq:chap:5p2.68}
       FNR(\theta,\psi) &= \frac{\mathbb{E}_{\theta}(\sum_{i:\theta_i \neq 0}(1-\psi_i))}{q_n} \nonumber \\
       &= \frac{1}{q_n} \sum_{i:\theta_i \neq 0} t^{EB}_{2i}. 
   \end{align}
   Recall that, in \eqref{eq:chap5p2.10}, we have already derived that, for each $|\theta_i| \geq \sqrt{2 \log (\frac{n}{q_n}) }+b$, $i \in S_{\boldsymbol{\theta}}$,
    \begin{align*}
        t^{EB}_{2i} 
       & \leq 1- \Phi(b)+o(1),
    \end{align*}
    where the $o(1)$ term is independent of any $i$. This combining with \eqref{eq:chap:5p2.68} implies for all sufficiently large $n$
    \begin{align} \label{eq:chap:5p2.69}
         \sup_{\theta \in \Theta_b} FNR(\theta,\psi) & \leq 1- \Phi(b)+o(1). 
    \end{align}
    Finally, \eqref{eq:chap:5p2.67} and \eqref{eq:chap:5p2.69} imply the upper bound to the left hand side of \eqref{eq:chap:5p2.62} is of the order $1- \Phi(b)+o(1)$ and completes the proof.
\end{proof}

\begin{proof}[Proof of Theorem \ref{chap5p2minmaxFDRFNRFB}]
   Here, we also follow the basic architecture of the proof of Theorem \ref{chap5p2minmaxFDRFNRtun}. Hence, in this case, too, it is sufficient to show that the upper bound to the left-hand side of \eqref{eq:chap:5p2.620} is of the order $1-\Phi(b)+o(1)$. \\
   Note that, using the argument of \eqref{eq:chap:5p2.52} and \eqref{eq:T-6.3a}, we have
 $$u_{7n}= \mathbb{E}_{\boldsymbol{\theta}} \bigg[\frac{\sum_{i:\theta_i=0}\psi_{i}^{\textrm{FB}}}{\max(\sum_{i=1}^{n}\psi_{i}^{\textrm{FB}},1)} \mathbf{1}\left\{\sum_{i=1}^{n}\psi_{i}^{\textrm{FB}} > \lambda q_n\right\}\bigg],$$ and
   $$u_{8n}= \mathbb{E}_{\boldsymbol{\theta}} \bigg[\frac{\sum_{i:\theta_i=0}\psi_{i}^{\textrm{FB}}}{\max(\sum_{i=1}^{n}\psi_{i}^{\textrm{FB}},1)} \mathbf{1}\left\{\sum_{i=1}^{n}\psi_{i}^{\textrm{FB}} \leq \lambda q_n\right\}\bigg],$$
   with $\psi_{i}^{\textrm{FB}}$ being defined in (23). 
   Note that
   \begin{align} \label{eq:chap:5p2.630}
        u_{7n} & \leq \frac{\sum_{i:\theta_i=0}t^{FB}_{1i}}{\lambda q_n} \nonumber \\
       & \lesssim \frac{n \alpha_n}{q_n \sqrt{\log (\frac{1}{\alpha_n}) }} (1+o(1)) \nonumber \\
        &= o(1), \textrm{ as } n \to \infty. 
    \end{align}

   Since the ratio term in the definition of $u_{8n}$ is bounded above by one,
   \[
   u_{8n}
   \leq
   \mathbb{P}_{\boldsymbol{\theta}}\left(\sum_{i=1}^{n}\psi_i^{\textrm{FB}}\leq \lambda q_n\right)
   \leq
   \mathbb{P}_{\boldsymbol{\theta}}\left(\sum_{i:\theta_i\neq 0}\psi_i^{\textrm{FB}}\leq \lambda q_n\right).
   \]
   Let
   \[
   \psi_i^{(0)}
   =
   \mathbf{1}\left\{
   \mathbb{E}(1-\kappa_i\mid X_i,1/n)>\frac12
   \right\}.
   \]
   Since
   \[
   \mathbb{E}(1-\kappa_i\mid X)
   =
   \int_{1/n}^{\alpha_n}
   \mathbb{E}(1-\kappa_i\mid X_i,\tau)\pi(\tau\mid X)\,d\tau
   \geq
   \mathbb{E}(1-\kappa_i\mid X_i,1/n),
   \]
   it follows that
   \[
   \psi_i^{\textrm{FB}}\geq \psi_i^{(0)}, \qquad i=1,\ldots,n.
   \]
   Hence,
   \[
   \left\{
   \sum_{i:\theta_i\neq0}\psi_i^{\textrm{FB}}\leq \lambda q_n
   \right\}
   \subseteq
   \left\{
   \sum_{i:\theta_i\neq0}\psi_i^{(0)}\leq \lambda q_n
   \right\}.
   \]
   Let
   \[
   S^{(0)}
   =
   \sum_{i:\theta_i\neq0}\psi_i^{(0)}
   \]
   and denote by $t^{(0)}_{2i}$ the type II error probability corresponding to the fixed rule $\psi_i^{(0)}$. Then
   \[
   \mathbb{E}_{\boldsymbol{\theta}}[S^{(0)}]
   =
   \sum_{i:\theta_i\neq0}(1-t^{(0)}_{2i}).
   \]
   The variables $\psi_i^{(0)}$, $i\in S_{\boldsymbol{\theta}}$, are independent Bernoulli random variables. Therefore, applying Hoeffding's inequality,
   \[
   u_{8n}
   \leq
   \mathbb{P}_{\boldsymbol{\theta}}\left(S^{(0)}\leq \lambda q_n\right)
   =
   \mathbb{P}_{\boldsymbol{\theta}}\left(S^{(0)}\leq \mathbb{E}_{\boldsymbol{\theta}}[S^{(0)}]-t\right),
   \]
   where
   \[
   t=q_n-\sum_{i:\theta_i\neq0}t^{(0)}_{2i}-\lambda q_n.
   \]
   Using the same type II error bound as in the proof of Theorem \ref{chap5p2minmaxFB}, applied to the fixed rule $\psi_i^{(0)}$, we obtain for sufficiently large $n$
   \begin{align}
        t &= q_n -\sum_{i:\theta_i \neq 0}t^{(0)}_{2i}-\lambda q_n \nonumber\\
     & \geq [\Phi(b)-\lambda+o(1)]q_n. 
   \end{align}
   This lower bound with Hoeffding's inequality confirms, for sufficiently large $n$
   \begin{align}\label{eq:chap:5p2.660}
       \sup_{\theta \in \Theta_b} u_{4n} =o(1), \textrm{ as } n \to \infty. 
   \end{align}
   Combining \eqref{eq:chap:5p2.630} and \eqref{eq:chap:5p2.660} ensures
   \begin{align}\label{eq:chap:5p2.670}
       \sup_{\theta \in \Theta_b} FDR(\theta,\psi) =o(1), \textrm{ as } n \to \infty. 
   \end{align}
   For $FNR(\theta,\psi)$, from definition, we obtain
   \begin{align}\label{eq:chap:5p2.680}
       FNR(\theta,\psi) &= \frac{\mathbb{E}_{\theta}(\sum_{i:\theta_i \neq 0}(1-\psi_i))}{q_n} \nonumber \\
       &= \frac{1}{q_n} \sum_{i:\theta_i \neq 0} t^{FB}_{2i}. 
   \end{align}
   Recall that, in Theorem \ref{chap5p2minmaxFB}, we have already derived that, for each $|\theta_i| \geq \sqrt{2 \log (\frac{n}{q_n}) }+b$, $i \in S_{\boldsymbol{\theta}}$,
    \begin{align}
        t^{FB}_{2i} 
       & \leq 1- \Phi(b)+o(1),
    \end{align}
    where the $o(1)$ term is independent of any $i$.
    This combining with \eqref{eq:chap:5p2.680} implies for all sufficiently large $n$
    \begin{align} \label{eq:chap:5p2.690}
         \sup_{\theta \in \Theta_b} FNR(\theta,\psi) & \leq 1- \Phi(b)+o(1). 
    \end{align}
    Finally, \eqref{eq:chap:5p2.670} and \eqref{eq:chap:5p2.690} imply the upper bound to the left hand side of \eqref{eq:chap:5p2.620} is of the order $1- \Phi(b)+o(1)$ and completes the proof.
\end{proof}

\end{document}